\documentclass[authoryear,preprint,review,12pt]{elsarticle}

\usepackage[toc,page]{appendix}
\usepackage{booktabs}
\usepackage{chngpage}
\usepackage{supertabular}
\usepackage{longtable,tabu}
\usepackage{hhline}
\usepackage{multirow}
\usepackage{multicol}
\usepackage{eurosym}
\usepackage{tabularx}
\usepackage{blindtext}
\usepackage{enumerate}
\usepackage{color}
\usepackage{array}
\usepackage{bm}
\usepackage[mathscr]{eucal}
\usepackage{mathtools}
\usepackage{amsfonts}
\usepackage{empheq}
\usepackage{nomencl}
\usepackage{colortbl}
\usepackage{setspace}
\usepackage[T1]{fontenc}
\usepackage[utf8]{inputenc}
\usepackage{geometry}
\usepackage[linesnumbered,ruled,vlined]{algorithm2e}
\usepackage{amsthm,amsmath,amssymb}
\usepackage{float}
\usepackage{subfig}
\usepackage{mathrsfs}
\usepackage{natbib}
\usepackage{graphicx}
\usepackage[colorlinks]{hyperref}

\graphicspath{{PICS1/}}

\begin{document}
 \title{Improving Critical Buildings Energy Resilience via Shared Autonomous Electric Vehicles - A Sequential Optimization Framework\tnoteref{t1,t2}}
\author{Jinming Liu\fnref{fn2}}

\author{Adam F. Abdin\corref{cor1}\fnref{fn1}}
\ead{adam.abdin@centralesupelec.fr}

\author{Jakob Puchinger\fnref{fn3,fn1}}

\cortext[cor1]{Corresponding author}

 \affiliation[fn2]{organization={Ecole Centrale Pekin, Beihang University}, 
                 addressline={37 Xueyuan Road},
                 postcode={100191}, 
                 city={Beijing}, 
                 country={China}}

 \affiliation[fn1]{organization={Universit\'e Paris-Saclay, CentraleSup\'elec, Laboratory of Industrial Engineering},
                 addressline={3 rue Joliot Curie}, 
                 city={Gif-sur-Yvette},
                 postcode={92290}, 
                 country={France}}

\affiliation[fn3]{organization={EM Normandie Business School, M\'etis Lab},
                 addressline={30-32 rue Henri Barbusse}, 
                 city={Clichy},
                 postcode={92110}, 
                 country={France}}
                 
 \begin{abstract}
The interdependence between electric power systems and transportation systems is rapidly increasing due to the high adoption of Electric Vehicles (EVs) and their charging infrastructures. Electric vehicles can represent additional load for the power system, but can also bring new opportunities for contributing to the efficient and resilient operations of the power grid. This is mainly because of their ability to provide back power to the system when it is not used for transportation, essentially serving as a moving battery source for the power grid. This so-called Vehicle-to-Grid (V2G), Vehicle-to-Building (V2B), or, generally, Vehicle-to-X (V2X) capability of EVs has been extensively studied in the literature. However, the upcoming development of autonomous driving systems and their integration within sharing mobility services can significantly add to the possibilities of interaction between the transportation system and the power grid. This paradigm is studied to a much lesser extent in the existing literature. Shared Autonomous Electric Vehicles (SAEVs) could allow for more control of the actions of the fleet allowing for large-scale coordinated responses both to mobility and energy demands. This coordinated response can be particularly useful for providing emergency power services in case of power loss while maintaining a high level of transportation service. Thus, improving the overall resilience of the system. In this work, we develop a dynamic optimization framework to evaluate the potential contribution of the SAEV fleet for improving critical buildings' energy resilience via V2B services. The model considers passengers' pick-up and transportation, relocation of vehicles, and battery charging and discharging. Power outage scenarios for critical buildings are considered and the potential of the SAEV fleet to fully or partially respond to the emergency power outage is studied. In addition, sensitivity analysis for key parameters such as the outage parameters is introduced. The results of the case study for the Ile-de-France region in France shows that the SAEV fleet has the potential to provide V2B service for critical building at an acceptable loss of passenger total waiting time on the transportation side. Furthermore, it is shown that it is capable of satisfying the emergency power load at a lower cost compared to investing in extra backup generators unless the outage occurs at significantly high frequencies.

 \end{abstract}

\maketitle

\section{Introduction}

Electric vehicles (EVs) are being rapidly adopted globally as worldwide policies are in favor of this change \citep{hertzke2017dynamics}. Their low operational costs, low greenhouse gas (GHG) emissions, and weak dependency on oil-based fuels make them a widely preferred choice for sustainable mobility systems, compared to traditional combustion vehicles \citep{chen2016management}. This led to more than 10 million electric cars being deployed on the roads worldwide in 2020 \citep{plotz2014will}. The increasing number of EVs and their necessity for electric power charging has strongly strengthened the interdependence between the electrified mobility system and the electric power system. This interdependence comes not only from more load being demanded from the power grid but also due to EVs' potential for providing ancillary grid services, such as Vehicle-to-Grid (V2G) power and its derived technologies.

Autonomous driving is another technology that is expected to have a high adoption rate in the near future \citep{rosenzweig2015review}. Level-2 automation (L2 level), which represents partial driving automation has already been implemented while level-3 automation, which represents conditional driving automation has been largely developed. For example, for highway assistant for Audi A8 and the Google cars \citep{brenner2018overview}. As for the fully driving automation (level-4 automation), it has been tested at low speeds on closed-off sections or testing grounds \citep{brenner2018overview}. Fully autonomous driving would lead to high environmental benefits due to the improvements in fuel economy, the reduction of the number of vehicles on the roads, and the reduction of the number of parking spaces needed as vehicles can operate for much longer hours without the need for a driver \citep{eugensson2013environmental}. This would also lead to high productivity gains for society, a decline in the number of accident and death tolls, and stress reduction \citep{rosenzweig2015review}.

Moreover, mobility systems are undergoing another important development related to the rise in the sharing economy. Shared mobility, with services such as Car2Go and ZipCar, is widely spreading across Europe and North America \citep{COPPOLA20191, chen2016carsharing}. Conventional vehicles are utilized for less than 10\% of their service time \citep{fagnant2014travel}. Conversely, shared cars have shown significant potential for enhanced convenience, energy efficiency, and economic advantages when contrasted with conventional traveling approaches. In several cases, shared cars have been shown to have the capability to replace around 9 traditional vehicles in providing an equivalent level of mobility service \citep{greenblatt2015automated}. With these developments, researchers and practitioners have recognized the potential synergy in implementing autonomous driving technology within shared mobility services \citep{narayanan2020shared}. The introduction of shared autonomous vehicles (SAVs) is, thus, expected to further improve the service efficiency and contribution towards effectively reducing overall GHG emissions \citep{fagnant2014travel}.

SAVs have significant potential to resolve some practical limitations of today’s traditional electrified mobility systems. This includes reduction of range anxiety, easier access to charging infrastructure, and efficient charging time management \citep{chen2016operations}. In addition, as more considerations are given to sustainability and environmental aspects of the mobility sector, SAVs are expected to, also, be largely driven by electric power, resulting in a significant deployment of shared autonomous electric vehicles (SAEVs). Electrically powered fleets of SAEVs have the potential to serve as ancillary providers to the power grid. This is achievable through their Vehicle-to-Grid (V2G) or Vehicle-to-Building (V2B) capabilities. Moreover, by adequately coordinating their recharging schedules in accordance with real-time charging costs, these fleets might also generate profits \citep{liao2021shared}. In addition, privately operated EVs are distributed resulting in high uncertainties in their V2G or V2B service patterns \citep{clement2009impact}, whereas a centrally operated SAEV fleet can ensure both mobility and power grid services homogeneously by its control pattern.

While several works have studied the economic benefits and service stability of SAEV fleets, further research on the potential impact of an interdependent SAEV transportation system and the power grid is required. This topic consists of a variety of situations, among which, is the possibility of SAEVs fleets to improve the power grid resilience by responding to emergency cases of power grid outages. Power grid outages are a rising phenomenon that causes a large number of adverse effects. For example, extreme weather events are one of the leading causes of power grid failures which have seen a significant rise in recent years \citep{kenward2014blackout, Abdin2019}. These outages do not only cause minor inconveniences in people's daily lives but also have a severe impact on critical infrastructures, causing severe impairment in delivering important services such as manufacturing of key products or the treatment of life-threatening events in hospitals. Therefore, it has become increasingly important to consider adequate solutions for mitigating and responding to these critical events while considering how to optimally plan for the desired cost and level of system resilience. 

On the electric power system side, backup generators have been the common solution for responding to emergency outage situations. However, these would typically require a high investment cost as well as high installation and operational costs, for a limited number of outage hours per year. Providing emergency power services directly to buildings via V2B capabilities could be a viable efficient alternative for improving the resilience of critical buildings against power outages during the outage duration. However, traditional V2G and V2B services that are provided by distributed and individually owned EVs are often limited by the fleet size, uncertainties in the charging behaviors of vehicle owners, and service stability. The large adoption of SAEV fleets, which are centrally operated, not subject to a human driver limitation, and continuously in service, can potentially ensure the right balance between providing critical grid services in case of emergency and meeting the mobility demands for which they are originally intended. However, to our knowledge, the study of this paradigm and the validation of its feasibility is still largely missing in the literature. 

Owing to the above considerations, the aim of this work is to develop an optimization modeling framework for the study of the potential of SAEV to improve the resilience of the interdependent power system - transportation system  by responding to emergency power outages in critical buildings via V2B power capabilities, while maintaining the highest level of service for mobility demands. To achieve this, an operational model of the SAEV fleet is developed, considering passenger pick-up and transportation, relocating, and battery charging/discharging. Moreover, power outage scenarios for critical buildings are modeled considering the energy demand and the infrastructure cost of backup sources. Then a V2B operation is introduced to model the potential of the SAEV fleet to fully or partially respond to the emergency power outage by relocating to the outage location and contributing to the power demand of the critical building. In this case, the model should be capable of optimally finding the best fleet dispatch to satisfy the emergency power demand considering the vehicles' state of charge (SOC), relocation time, battery consumption, duration of intervention, and the satisfaction of passenger transportation requirements, among other realistic considerations.

Towards this goal, we propose a novel optimization model for SAEV fleet operations considering dynamic travel demands, charging requirements of SAEV batteries, and V2B capabilities for responding to critical building power outages in cases of emergency. The model considers the main objective of minimizing the passengers total waiting time as the primary goal of the SAEV fleet. To ensure operational efficiency, penalty terms are added to the objective that ensures the minimization of the relocation time of vehicles and the total charging costs for the SAEVs batteries. We assume a fleet of vehicles fully owned by one operator for providing an Automated Mobility-on-Demand (AMoD) service. A realistic case study based on the travel demand patterns of the Ile-de-France region (larger Paris area in France) is considered for the application of the proposed approach and the analysis of the results.  First, we consider a SAEV operational scenario under normal conditions (i.e. no emergency power outages in critical buildings) as the base case for the analysis. In this scenario, a dynamic vehicle routing problem is solved to obtain the optimal SAEV plan to relocate to pick up passengers and transport them to their respective destinations. Battery charging decisions are, also, taken according to the energy requirements of the vehicles. Then, power outage scenarios are introduced for critical buildings, particularly, large hospitals. In this case, in addition to meeting travel demands and minimizing travelers waiting times, the SAEV fleet can use their available battery power to provide temporary power services to the critical building. We assume a power outage to occur for a limited number of hours in a district and we focus on one critical building (a large hospital). The SAEV fleet is, thus, regarded as a backup power source, and the vehicles are partially or totally required to relocate to the building to provide emergency power services while ensuring passengers traveling demand and transportation constraints.

\noindent The main contributions of this paper can be summarized as follows:
\begin{itemize}
  \item This paper proposes a novel optimization framework for the analysis of the interdependent power system resilience and SAEV mobility system. While considering the impact of SAEV on power grid operations is not a new topic, most of the related research focuses on normal operating conditions and on providing grid services such as load peak-shaving and frequency regulation. However, the potential of SAEV to provide emergency power services, which can have various important consequences on grid investments and operations, is still largely unexplored.  
  \item This paper proposes a MILP formulation for the dynamic Autonomous-Mobility-on Demand problem considering Vehicle-to-Building (V2B) power exchange capabilities for an SAEV fleet. Compared to existing models, novel extensions are proposed to account for power outages in critical buildings and the requirement for the SAEV fleet to provide backup power services while maintaining high transportation service levels. Furthermore, a Model Predictive Control (MPC) strategy is proposed to solve the dynamic optimization problem efficiently.
  \item A realistic case study based on travel data for the Ile-de-France region is considered to provide insights into the usefulness and applicability of the proposed approach. Historical mobility demand data are collected for the Ile-de-France region \citep{horl2021open} and are used to generate the transport demand distribution. Moreover, the average daily energy consumption of a hospital in Paris is used to build the power outage emergency model. The case study enables the comparison between investing in traditional backup generators as opposed to relying on a mixed generator-V2B backup power for improving the power system resilience. This case study can be used as a benchmark in the literature for further analysis of the interdependency between power systems and electrified mobility systems.

\end{itemize}

This article is structured as follows. In Section 2, a literature review is presented, covering the literature on SAEV mobility systems and the resilience evaluation for energy systems. In Section 3 the model proposed is presented in detail. This includes the proposed MILP model for the transportation layer, the MPC control algorithm used to solve the dynamic optimization problem, and an extension of the model to consider the case of a power outage. Section 4 presents the case study, including transport data, vehicle parameters, and power outage assumptions. In section 5, the optimization results are presented and discussed. Finally, Section 6 provides the conclusions and managerial recommendations.

\section{Literature review}

In this review, the literature on the methods to improve the resilience of power grids and provide backup power during emergencies is thoroughly analyzed. Then, the attention shifts to studies centered on V2B services in a general context. Lastly, the research on SAEV mobility systems is summarized, and separated into two sub-sections. The literature on SAEV fleet control and infrastructure planning is discussed in Section \ref{subsection-saevcontrol}, while the interdependence between SAEV systems and the power grid is explored through a review of relevant research works in Section \ref{subsection-saevinterdependency}.

\subsection{Power system resilience}
Several studies in the literature have explored the assessment of resilience in electric power systems. For instance, the study in \citep{panteli2016power} focuses on modeling the failure of electric networks during strong wind conditions and estimating the impact of enhancing the power grid infrastructure by adding transmission lines. Similarly, the work in \citep{cadini2017modeling} deals with the estimation of power grid resilience through the discussion of power line repair efficiency. The study in \citep{zhu2021spatio} utilizes real power outage data and proposes a spatiotemporal random process to analyze power grid resilience by considering identified resilience capabilities and interdependence between service areas. The study in \citep{liu2016microgrids} presents a framework for analyzing the resilience of an electric power grid with integrated microgrids during extreme conditions. The mesh grid approach is adopted and the impact of extreme events on power grid operations, and the relationship among the transmission grid, distribution grid, and microgrid is analyzed. The work of \citet{movahednia2021power} proposes a stochastic resource allocation approach for protecting power substations against flood events a day in advance. The research in \citep{abdin2019modeling} presents a modeling and optimization framework for long-term power systems planning aimed at ensuring the resilience of the power system against extreme heat waves and drought events. The study in \citep{alkhaleel2022model} proposes a two-stage mean-risk stochastic restoration model, using a MILP model, to study the post-disruption restoration problem for interdependent critical infrastructure networks under uncertainty, with the goal of minimizing the total cost associated with unsatisfied demands, repair tasks, and flow of interdependent infrastructure networks. \citet{widiputra2020optimal} present an algorithm to determine the optimal reconfiguration of the power system to increase grid resilience, ensuring that the operational parts of the power system remain functional during extreme weather events. \citet{hajializadeh2021rv} develop a framework for a Resilience and Vulnerability-Informed Decision Support System (RV-DSS) to identify the vulnerable components of interdependent networks, proposing a measure of network resilience in response to hazardous events.

On the supply side, \citet{chen2022demand} consider the improvement of power system resilience by studying the demand response  of household power consumption for a case study in China. The study also introduces mechanisms for allocating national targets for demand response development. In addition, the study proposes a framework designed to estimate the resulting economic benefits. \citet{wei2016quantifying} assess the impact of road accessibility loss on electric transmission and distribution networks, establishing two vulnerability analyses for both the transmission network and the distribution grid. The study concludes that the interdependence of transportation and electricity infrastructures must be considered in the design of the electricity network, though it mainly focuses on the design and identification of critical components of the power grid and does not address the potential for EVs to provide power services to the grid in emergency situations.

All of these studies, however, do not consider the potential interaction between the power system and electrified mobility systems and the potential improvements in resilience as a result of properly operated interdependent systems.

\subsection{Vehicle-to-Building (V2B) services }
A significant number of studies have investigated the potential of electrified mobility systems to provide Vehicle-to-Grid (V2G) services (see for example the works by \citet{ota2011autonomous, li2015optimizing, quddus2019modeling, iacobucci2019optimization, rajamand2020vehicle, bhatt2020review, das2020charging, yu2022electric, demirci2023comprehensive, alfaverh2023optimal, aktar2023framework}). However, the potential of EVs to directly provide power to buildings has received limited attention. This section focuses on the literature that has considered Vehicle-to-Building (V2B) services.

The use of EVs as a means to provide V2B services has been explored in the literature as discussed in detail in this section. For instance, a linear programming approach to generate the optimal charging/discharging schedule for an EV plugged into a smart building is presented in \citep{molina2012optimal}, under various time-varying electricity pricing schemes. The potential economic benefits of this method are demonstrated, particularly with the integration of energy storage systems.

In \citep{barone2019building}, a novel energy management system for buildings connected to a microgrid is proposed, with EVs considered as an active component of the microgrid's energy system. The study balances both power and mobility needs and shows that the proposed system improves building grid reliance and reduces grid electricity consumption. A scheduling optimization framework to reduce the peak load of an institutional building within a microgrid is presented in \citep{ouammi2021peak}. V2B capabilities of EVs are considered to control the operation of the microgrid and manage power flows. The results demonstrate that EVs can significantly contribute to peak load reduction and microgrid operational flexibility.

The topic of reducing peak load via V2B services has also received attention. In \citep{pang2011bevs}, EVs are used as energy storage to support the power grid, with the benefits and implementation issues of using EVs for demand side management and outage management discussed. \citet{odkhuu2018optimal} present an optimal energy management algorithm to minimize peak load on a university campus with V2B, coordinating EV charging/discharging to reduce peak load and considering real-time pricing. \citet{Geojson} propose a control algorithm to manage a commercial building's peak demand and daily electricity cost via EV V2B management.

The impact of EVs on grid resilience is considered in \citep{rahimi2018electric}. A Vehicle-to-Home (V2H) service is introduced in the case of power outages, and a min-max linear programming model is built to minimize the maximum load in different household electricity usage scenarios. The results show that EVs can provide a valuable solution in providing power during outages, but only as backup batteries and not considering transportation.

More recently, \citet{chai2023two} developed a two-stage optimization framework to determine the charging and discharging schedule for EVs participating in a V2G program at an office building. The approach developed seeks to minimize the total operating cost of the office building from the building owner’s perspective. The initial step involves employing day-ahead optimization. This process relies on the anticipated building load pattern and EV behavior. This serves the purpose of economizing on electricity expenses by minimizing the building's peak demand. Subsequently,  real-time optimization is executed to fine-tune the V2G operations. 

While previous studies discuss the potential for EVs to provide V2B services, the coordinated potential of SAEVs fleets in providing resilience power services while maintaining transportation service levels is not addressed. SAEVs offer a unique opportunity for electricity services in cases of power outages, due to their size and scope of operations, and different mobility considerations necessitate specific modeling frameworks and operational considerations.

\subsection{Scheduling shared EVs and SAEVs: fleet control and infrastructure design}
\label{subsection-saevcontrol}

The literature has explored various models for optimizing the operations of autonomous electric vehicle (SAEV) fleets. \citet{zhang2016model} propose a dynamic Autonomous-Mobility-on-Demand (AMoD) model that takes into account charging decisions and traffic demands to meet mobility requirements in a given area. The model is expressed as a mixed-integer linear programming (MILP) problem and solved using a Model Predictive Control (MPC) algorithm. \citet{lee2020shared} present a SAEV model that combines fleet control and charging station design, formulated as a sequential quadratic programming (SQP) problem. The study considers uncertainties and reliability-based design optimization (RBDO) principles. \citet{wang2022resilience} propose a mathematical framework to analyze the resilience of an electrified road network (ERN) under potential failures of its fast-charging stations (FCSs) and to assess the impact of FCS deployment locations on the charging infrastructure resilience.

In \citep{fagnant2014travel}, an agent-based SAEV model is developed, considering simplified traffic assumptions and a gridded road network. The study shows that a SAEV system can provide sufficient mobility services with a fleet size of about one-tenth the number of private vehicles required for equivalent quality of service. \citet{ma2021vehicle} propose a MILP model to optimize the service routes of an SAEV fleet and minimize the weighted objective of total travel distance, travel time, and energy consumption for the vehicles, taking into account charging and repositioning. \citet{folkestad2020optimal} consider charging and repositioning for the optimal assignment of a SAEV fleet. \citet{sheppard2019joint} propose a quadratic programming model to optimize the size of an SAEV fleet, charging stations, and their locations, with simulations conducted using relevant US data. The simulation result showed that 1 SAEV could replace around 10 personally owned vehicles with acceptable energy demand and peak charging load.

Moreover, \citet{melendez2020optimal} explore the optimal operation of a large SAEV fleet and a set of grid-supported charging centers throughout the transportation network, formulated as a robust MILP model. The objective is to maximize the fleet's profit, taking into account practical considerations, such as day-ahead load commitment and energy arbitrage. A case study demonstrates the ability to make operational decisions one day in advance and control SAEV dispatch and charging decisions in real time. \citet{zhang2022heterogeneous} present a MILP formulation for an EV car-sharing system that maximizes the operator's profit and demonstrates the advantage of considering a heterogeneous fleet and partial recharging for system flexibility, based on a case study in Beijing. \citet{lu2021optimal} study the optimal allocation of a fleet of plug-in EVs to the stations of an EV-sharing system, considering uncertain demand and maximizing the system operator's profit.

Despite these studies offering various methods for modeling and optimizing SAEV operations and infrastructure design, they do not take into account the interdependence between mobility systems and the power system, nor the potential of SAEVs in providing emergency power services for buildings via V2B services. Similar to the work proposed in \cite{zhang2016model}, our proposed sequential optimization model is formulated as a MILP and solved using MPC. However, our proposed modeling framework extends the formulation to consider V2B decisions (discharging capabilities in specific nodes), as well as novel constraints that describe the behavior of the interdependent systems during a power outage. Additionally, compared to the above-mentioned literature, our framework proposes novel metrics to evaluate the economic value of providing emergency power service via V2B, compared to traditional backup generators.

\subsection{SAEV: interdependence between the SAEV system and the power grid}
\label{subsection-saevinterdependency}

The interdependence between the electrified transportation system, particularly shared autonomous electric vehicle (SAEV) services, and the power grid has garnered attention in the literature. In particular, studies have aimed to understand the relationship between the two systems through various modeling and optimization frameworks.

For instance, \citet{fernandez2021inverse} introduce a data-driven two-step estimation approach that utilizes inverse optimization and kernel regression to determine the non-linear relationship between an aggregate price-response and its explanatory variables, thereby deriving a bid/offer curve for the charging of an EV fleet. \citet{zhang2020smart} propose a smart charging framework that aims to shift electricity demand from high-price peak periods to low-price renewable generation periods. The study employs an agent-based model and finds that larger battery capacities are more responsive to low electricity cost charging opportunities and have the potential to reduce energy-related costs.

\citet{iacobucci2018modeling} present a model combining SAEV fleet assignment and charging station design, which considers the SAEV fleet as a backup for the power grid operation via V2G service. The results of the study, which is based on a case study of Tokyo, show the influence of parameters such as the number of SAEVs and the location of charging stations on system performance, as well as the impact of V2G service on passengers' total waiting time and rejection rates.

\citet{iacobucci2019cascaded} introduce an optimal SAEV dynamic traffic assignment model that takes into account V2G services capabilities, modeled as a MILP and solved using a MPC algorithm. \citet{iacobucci2019optimization} extend the model to consider aggregate battery capacity and separate the transport operation and charging operation into two separate layers. Both studies examine the feasibility of providing V2G services by an SAEV fleet. 

More recently, studies started to explore the utilization of model-free decision support methods, particularly Reinforcement Learning (RL) to address the complexities of managing EV charging in uncertain conditions \citep{abdullah2021reinforcement}. For example, \citet{mhaisen2020real} investigated the use of RL to control the charging and discharging operation of EVs in a smart grid that supports V2G schemes considering real-time pricing of electric power. In \citep{hao2023v2g} a DQN-based RL method is proposed to address the charging management of EVs considering the power grid and the possible risks of overburdening the power system due to uncontrolled charging. A similar problem has been, also, investigated in \citep{alfaverh2023optimal} using a Deep-RL approach to find an optimal V2G control strategy for Supplementary Frequency Regulation (SFR) of the power system. 

While these studies have explored the interdependence between the SAEV system and the power grid in terms of coordinated charging scheduling, the potential of SAEVs to improve the resilience of the power demand in emergency situations has not been  discussed. In terms of modeling and solution algorithms, while the framework proposed in our study closely follows the work proposed in \citep{zhang2016model} and extended in \citep{iacobucci2019optimization}, our study is focused on modeling V2B services and the dynamic capability of the SAEV fleet to provide electricity to critical buildings during a power outage. To the best of our knowledge, this problem has not been treated in the existing literature. We also show that, despite the modeling differences, MPC remains a suitable algorithm that can be readily adjusted to handle the dynamic nature of these decisions with adequate computational time. MPC relies on a well-understood model of the system under investigation, and when accurate models are available, MPC can prove to be highly efficient. 

In light of this, this work proposes a dynamic AMoD optimization model that aims to minimize passengers' total waiting time taking into account V2B requirements in emergency power outage circumstances, and a solution method based on MPC. To the best of our knowledge, this problem has not been studied in the literature, and the framework proposed in this work could provide a roadmap for future research on the potential for improving the resilience of interdependent infrastructure systems.

The proposed modeling framework, solution method, and the analysis of the results of a practical case study are presented in the next sections. 

\section{Methodology}

\subsection{Model formulation}
This section describes the proposed optimization model for managing an SAEV fleet to satisfy transportation demand while taking into account charging considerations. Furthermore, the V2B capabilities of the SAEV fleet are considered for providing electricity to critical buildings during a power outage. The proposed model is formulated as a MILP following the work by \citet{iacobucci2019optimization}. 

In our model, a discrete-time system with $N$ nodes and $K$ vehicles is considered and it is simulated in $T$ time steps. Let $\mathcal{N}$ represent the set of nodes ($|\mathcal{N}|=N$), $\mathcal{K}$ represent the set of vehicles ($|\mathcal{K}|=K$), and $\mathcal{T}$ denote the set of time steps ($|\mathcal{T}|=T$).
Table \ref{table:1} summarizes all the variables and parameters of the proposed model. The model is described by a set of equality and inequality constraints that describe the operation of the SAEV fleet, and an objective function.

{\footnotesize
\begin{longtable}{p{1.5cm}|p{8cm}}

\hline
\multicolumn{2}{c}{A. Sets and Indexes} \\\hline
 $i$, $j$ & Index of nodes ($\in \mathcal{N}$). \\
 $n$ & Index of node where a power outage occurs ($\in \mathcal{N}$). \\
 $k$ & Index of vehicles ($\in \mathcal{K}$). \\
 $t$ & Index of time-period $\left( \in \mathcal{T} \right)$. \\
 $\theta_i$ & Index representing the number of time-steps remaining until a vehicle reaches node $i$ ($\in \{0, \text{max}_j \{t_{ji}\}-1$\}).\\
 $\bar{\Theta}_i$ & Maximum distance to travel to node $i$ from all connected nodes $j$ ($=\text{max}_j t_{ji}-1$) .\\\hline

\multicolumn{2}{c}{ B. Parameters and Exogenous variables }   \\\hline
 $\mathcal{B}$ &The battery capacity of the vehicle.\\
 $\mathcal{C}_v$ & The infrastructure cost per vehicle.\\
 $\mathcal{C}_{v2b}$& Total V2B cost.\\
 $\mathcal{C}_i$ & V2B infrastructure cost per year.\\
 $\mathcal{C}_r$ & Relocation cost per power outage.\\
 $\mathcal{C}_e$ & Energy cost per power outage.\\
 $P_{ijt}$ & Number of passengers arriving at node $i$ with destination $j$ at time $t$.  \\
 $\sigma_{t}$ & Electricity price at time $t$  [euro/kWh].\\
 $outage_{it}$ & Binary parameter, $outage_{it}=1$ if power outage occurs at node $i$ at time period $t$ and $0$ otherwise.\\
 $f_{out}$ & Annual frequency of power outage in critical buildings.\\
 $Q_{m}$ & Energy generated by backup generators during a power outage.\\
 $Q_{d}$ & Energy required for the critical building during a power outage.\\
 $\gamma^{max}$ & Maximum State-of-Charge (SOC) of a vehicle [in \%].\\
 $\gamma^{min}$ & Minimum State-of-Charge (SOC) of a vehicle [in \%].\\
 $\gamma^{init}$ & Initial State-of-Charge (SOC) of a vehicle [in \%]. \\
 $t_{ji}$ & Number of time-steps for a vehicle to move from node $j$ to node $i$.  \\
 $\vartheta_{d}$ & Energy consumption rate of a vehicle while moving.\\
 $\vartheta_{c}$ & Charging rate of vehicles battery.\\
 $\vartheta_{v2b}$ & Discharging rate of vehicles battery. \\
 $\eta$ & Efficiency of V2B.\\
 $\rho_{1}$ & Penalty weight of rebalancing cost. \\
 $\rho_{2}$ & Penalty weight of charging cost. \\
 $\tau$ & Duration of each time-step $t$ of the model.\\
 $\omega$ & Cost of cycling the battery, incurred when supplying power to the grid.\\\hline

 \multicolumn{2}{c}{C. State Variables }    \\\hline 
 $d_{ijt}$/$D_{ijl}$ & Number of passengers waiting at node $i$ with destination $j$ in prediction time-step $(t)$ / real time-step $(l)$.  \\
 $ a^{k\theta_i}_{it}$/$ A^{k\theta_i}_{il}$ & Binary variable which defines the movement of vehicle $k$ in prediction time-step $(t)$ / real time-step $(l)$. \\
 $\gamma^k_{t}$/$\Gamma^k_{l}$ & The state of charge (SOC) of vehicle $k$ at prediction time-step $(t)$ / real time-step $(l)$.  \\
 $u^k_{it}$/$U^k_{il}$ & Binary variable representing the parking status of a vehicle $k$ at node $i$ in prediction time-step $(t)$ /real time-step $(l)$.  \\\hline
\multicolumn{2}{c}{D. Control Variables }    \\\hline
$e^k_{t}$ & Energy charged by vehicle $k$ at time $t$.  \\
$g^k_{t}$ & Energy discharged by vehicle $k$ at time $t$.  \\
$q_t$ & Energy received by the building at time $t$.  \\
$q_{v2b}$ & Total energy discharged by the SAEV fleet to the building ($= \mathcal{B} * \sum_{k,t} g_t^k$).\\
$v^k_{ijt}$ &Binary variable which is equal to 1 if vehicle $k$ is moving a passenger from node $i$ to node $j$ at time $t$, and 0 otherwise.\\
$r^k_{ijt}$ &Binary variable which is equal to 1 if vehicle $k$ is relocating from node $i$ to node $j$ at time $t$, and 0 otherwise.\\\hline
\caption{Model variables, parameters and constants.}
\label{table:1}
\end{longtable}}

First, the number of passengers $d$ waiting at node $i$ with the destination to node $j$ at time $t+1$ is defined to be equal to the number of passengers waiting in the previous time period, plus newly arrived passengers, minus the number of passengers that are picked up by a vehicle. This can be formulated as :
\begin{equation}
d_{ijt+1} = d_{ijt} + P_{ijt}- \sum\limits^K_{k=1}v^k_{ijt},{\forall}{i,j} {\in}{\mathcal{N}},t {\in}{\mathcal{T}}
\end{equation}

Since we are modeling dynamic mobility-on-demand service, it is important to keep track of the vehicle location within the system. This is done by keeping track of how long a vehicle will be traveling before it reaches its destination, as proposed in \citep{zhang2016model}. The evolution of vehicle movements is captured by a binary variable $a^{k\theta_i}_{it}$ which tracks the distance of the vehicle to node $i$ from its destination. $\theta_i \in \{0, \text{max}_j\{t_{ji}\}\}$ is the number of time steps remaining until the vehicle reaches node $i$ from neighboring node $j$.  The case that a vehicle arrives at node $i$ is captured when $a^{k0}_{it}=1$. The evolution of the vehicle movement $ a^{k\theta_i}_{it}$ is, then, described in two situations:
\begin{itemize}

\item If $\theta_i<\bar{\Theta}_i$, where $\bar{\Theta}_i$ is the maximum distance to node $i$, the vehicle location $a$ for vehicle $k$ from node $i$ with distance $\theta_i$ is obtained as the vehicle location in the previous time period (at a distance of $\theta_i+1$ from node $i$) if the vehicle was already on-route, or the departure of vehicle $k$ from node $j$ to node $i$ after a decision to carry a passenger or relocating, with the distance between these nodes being $\theta_i+1$. 

\item If $\theta_i=\bar{\Theta}_i$, $ a^{k\theta_i}_{it}$ is only obtained as the decision to pick up passengers or to relocate.
\end{itemize}

\noindent Thus, the evolution of the movement is as expressed in the following equation:
\begin{equation}
a^{k\theta_i}_{it+1}=
\left\{
             \begin{array}{lr}
             a^{k\theta_i+1}_{it}+\sum\limits_{j:t_{ji}-1=\theta_i}(v^k_{ijt}+r^k_{ijt}), &  \theta_i<\bar{\Theta}_i\\
             \sum\limits_{j:t_{ji}-1=\theta_i}(v^k_{ijt}+r^k_{ijt}), &  \theta_i=\bar{\Theta}_i
             \end{array}
\right. ,{\forall}{i}{\in}{\mathcal{N}},k{\in}{\mathcal{K}},t{\in}{\mathcal{T}},\theta_i{\in}{[0,\bar{\Theta}_i]}
\end{equation}

\noindent The evolution of the parking status $u$ of each vehicle $k$ at node $i$ and time $t+1$ is equal to the parking state of the vehicle at the previous time period, plus if the vehicle has just arrived $(a^{k0}_{it}=1)$, minus the case that vehicle leaves the node either after picking-up a passenger $(v^k_{ijt}=1)$ or for rebalancing $(r^k_{ijt}=1)$. This constraint is formulated as described in the following equation:
\begin{equation}
u^k_{it+1} = u^k_{it} + a^{k0}_{it}- \sum\limits_{j}(v^k_{ijt}+r^k_{ijt}),{\forall}{i}{\in}{\mathcal{N}},k{\in}{\mathcal{K}},t{\in}{\mathcal{T}}
\end{equation}

\noindent Similarly, the evolution of the charging state $\gamma$ of the battery of each vehicle $k$ at time $t + 1$ is equal to the previous charging state of the vehicle at time $t$, plus the energy charged $e_t^k$, minus the energy discharged $g_t^k$	 minus the energy consumed while the vehicle is moving $\vartheta_d \sum\limits_i^{}\sum\limits^{\bar{\Theta}_i}_{\theta_i=0}  a^{k\theta_i}_{it+1}$. This is given by the following equation:
\begin{equation}
	\gamma^k_{t+1} = \gamma^k_{t} + e^k_t- g^k_t-\vartheta_d \sum\limits_i^{}\sum\limits^{\bar{\Theta}_i}_{\theta_i=0}  a^{k\theta_i}_{it+1},{\forall} k{\in}{\mathcal{K}},t{\in}{\mathcal{T}} 
\end{equation}

\noindent Notice that the vehicle cannot charge and discharge at the same time since an optimal solution will require that either charging or discharging occurs concurrently.

\noindent Next, logical constraints are specified. First, for each time step $t$, a vehicle $k$ can either be parked at one node ($\sum\limits_i^{}u^k_{it}=1$) or moving towards one node ($\sum\limits_i^{}\sum\limits^{\bar{\Theta}_i}_{\theta_i=0}  {}a^{k\theta_i}_{it}=1$) but not both. This is described in the following constraint:
\begin{equation}
	\sum\limits_i^{}u^k_{it}+\sum\limits_i^{}\sum\limits^{\bar{\Theta}_i}_{\theta_i=0}  {}a^{k\theta_i}_{it}=1,{\forall} k{\in}{\mathcal{K}},t{\in}{\mathcal{T}}
\end{equation}

\noindent Vehicles can only perform one task at each time period $t$. Vehicles can either be relocating or picking up a passenger. Moreover, a vehicle cannot be parked in the time period following the decision to pick up a passenger or relocate to another node. This relation is shown in the following constraint:
\begin{equation}
	\sum\limits_i^{}\left(u^k_{it+1}+\sum\limits_j\left(v^k_{ijt}\right)+\sum\limits_j\left(r^k_{ijt}\right)\right)\leq1,{\forall} k{\in}{\mathcal{K}},t{\in}{\mathcal{T}}
\end{equation}

\noindent Moreover, the number of vehicles that are requested to transport passengers is less than the number of passengers waiting plus the newly arriving passengers, as shown in the following constraint:
\begin{equation}
	\sum\limits_k v^k_{ijt}\leq d_{ijt}+P_{ijt},{\forall}{i,j}{\in}{\mathcal{N}}, t{\in}{\mathcal{T}}
\end{equation}

\noindent Regarding the vehicles’ electric energy status. Constraint (\ref{eq8}) limits the energy of each vehicle $k$ at time $t$ to be smaller than the maximum battery capacity. Moreover, available energy in the vehicle should be higher than a minimum state-of-charge ($\gamma^{min}$) plus the energy required for moving passengers ($\vartheta_{d}t_{ij}v^k_{ijt}$) or relocating the vehicle ($\vartheta_{d}t_{ij}\omega^k_{ijt}$) as in constraint (\ref{eq9}). This ensures that vehicles would not start a journey from node $i$ to node $j$ if the battery level is not enough to ensure this journey.
\begin{equation}
	\gamma^k_t\leq \gamma^{max},{\forall}k{\in}{\mathcal{K}}, {t}{\in}{\mathcal{T}}
\label{eq8}
\end{equation}

\begin{equation}
	\gamma^k_t\geq \vartheta_{d}t_{ij}v^k_{ijt}+\vartheta_{d}t_{ij}r^k_{ijt}+\gamma^{min},{\forall}{i,j}{\in}{\mathcal{N}},k{\in}{\mathcal{K}},t{\in}{\mathcal{T}}
\label{eq9}
\end{equation}

\noindent Charging and discharging constraints are expressed in Eq.~(\ref{eq10}) and (\ref{eq11}), respectively. Eq.~(\ref{eq10}) limits the charging rate to a maximum rate and ensures that vehicles can only charge when parked. Similarly, Eq.~(\ref{eq11}) limits the discharging rate to the maximum V2B rate and ensures that vehicles cannot discharge except when parked at a discharging node. Moreover, if a power outage occurs at a given node, vehicles cannot charge their batteries at this node, however, it is still possible that they provide power back to the grid or connected buildings. These two constraints are expressed in the following equations.
\begin{equation}
	e^k_t\leq \vartheta_{c}\sum\limits_i\left(\left(1-outage_{it}\right)u^k_{it}\right),{\forall} k{\in}{\mathcal{K}},t{\in}{\mathcal{T}}
\label{eq10}
\end{equation}

\begin{equation}
	g^k_t\leq \vartheta_{v2b}\sum\limits_i\left(outage_{it}u^k_{it}\right),{\forall} k{\in}{\mathcal{K}},t{\in}{\mathcal{T}}
\label{eq11}
\end{equation}

\noindent Next, the relation between the energy received by the critical building  ($q_t$) and the V2B energy discharge by each vehicle  ($g_t^k$) is given in Eq.~(\ref{eq12}). The energy received by the building at time $t$ is the sum of all the energy discharged by the vehicles at $t$ multiplied by the discharging efficiency $\eta$:
\begin{equation}
\label{eq12}
	q_t=\eta \sum\limits_k g_t^k,{\forall} t{\in}{\mathcal{T}}
\end{equation}

\noindent The proposed model in this paper considers the case of emergency power outages in critical buildings. In this case, the SAEV fleet could provide temporary backup electricity for a limited amount of time through the V2B capability. In case of a power outage, the energy obtained by the critical building via the SAEV fleet is ensured by Eq.(\ref{eq14}), as follows:

\begin{equation}
	q_t\geq outage_{nt}\left(Q_d-Q_m\right),{\forall}{t}{\in}{\mathcal{T}}
\label{eq14}
\end{equation}

\noindent Finally, all the variables should be non-negative.

\noindent \textbf{Objective Function:}

The objective of this model is to ensure the best possible mobility service level while contributing to improving the energy resilience of critical buildings in case of power outages.  The objective function is, therefore, aimed at minimizing the passengers' total waiting time $\xi^x$ (Eq.~\ref{eq15}) for the mobility demand services. In addition, two penalty terms: $\xi^u$, and $\xi^m$ are added to the objective function to ensure the most efficient operation of the fleet. The first term ($\xi^u$) represents the relocating time of empty vehicles (Eq.~\ref{eq16}) and the second term ($\xi^m$) represents the charging cost of the vehicles’ batteries (Eq.~\ref{eq17}). Those three terms are formulated as follows:
\begin{equation}
	\xi^x_t=\sum\limits_i\sum\limits_j d_{ijt},{\forall} t{\in}{\mathcal{T}}
\label{eq15}
\end{equation}
\begin{equation}
	\xi^u_t=\sum\limits_k\sum\limits_i\sum\limits_j t_{ij}r^k_{ijt},{\forall} t{\in}{\mathcal{T}}
\label{eq16}
\end{equation}
\begin{equation}
	\xi^m_t=\sum\limits_k \left( \left(e_t^k-\eta g_t^k \right)\sigma_t+\omega g_t^k\right),{\forall} t{\in}{\mathcal{T}}
\label{eq17}
\end{equation}

\noindent The full objective function, therefore, is given by Eq.~\ref{eq18}. Note that the primary objectives of incorporating the terms $\rho_1\xi_t^u$ and $\rho_2\xi_t^m$ are twofold: to prevent unnecessary vehicle rebalancing and to avoid excessive charging costs. In the absence of rebalancing limitations, vehicles would always be occupied by random rebalancing demands, rendering them unable to promptly respond to passenger requests. Additionally, in the absence of charging constraints, vehicles would charge whenever they are unoccupied, incurring an excessively high charging cost that is deemed unnecessary. Penalty terms are also employed to ensure that minimizing passenger waiting time takes precedence over the reduction of relocation time and charging costs. Nevertheless, the weight of the penalty terms is of minor significance, as its primary role is to guarantee the prioritization of minimizing passenger waiting time qualitatively. Thus, its weight should be substantially smaller than that of the main objective, while being within a reasonable order of magnitude of the other parameter values to ensure the numerical stability of the algorithm. Accordingly, to satisfy these criteria, we investigated a wide range of penalty value combinations covering the range from 0.01 to 1 for ($\rho_1$) with 0.01 step increase and the range from 0.001 to 0.01 for ($\rho_2$) with a 0.001 step increase and evaluated the impact of these penalties on the results. Our analysis showed that we can set the weight of the relocation cost ($\rho_1$) to 0.01, which is sufficiently low to have no significant effect on the major objective while achieving minimal relocation time in our model. Moreover, the weight of the charging cost ($\rho_2$) should be any value lower than $\rho_1$ to ensure necessary relocation and is thus set to 0.001.
\begin{equation}
	\min \sum\limits_{t=0}^{T} 	\left( \xi^x_t+\rho_1 \xi^u_t+ \rho_2 \xi^m_t\right)
\label{eq18}
\end{equation}
\subsection{Solution algorithm: Model Predictive Control (MPC)}
\label{solution_section}
MPC solves an open-loop optimization problem at each time step to provide a series of control actions in a ﬁxed horizon, and the ﬁrst control action is executed in each iteration. It is shown in \citep{zhang2016model} that, for autonomous mobility-on-demand (AMoD) routing and scheduling applications, MPC can achieve closed-loop performance, and can be run in real-time for moderately sized systems that can be encountered in practice. MPC method applies a receding horizon control that can account for complex decision-making in practical real-world scenarios. In the previous section, a multi-stage model was proposed with the aim of tracking the dynamic characteristics of the system. In order to solve this model, the MPC method was selected as a solution approach due to its ability to effectively handle multi-stage optimization problems.

To solve the proposed model within an MPC framework, first, we solve the MILP model described above at real-time instant $l$ for the given optimization horizon $T$. Next, the resulting state variables at the second time step are fixed as the initial conditions at real-time instance $l+1$ and we consider the control variables at the first-time step as the real control action at real-time instance $l$. For example, for a prediction horizon between $T=1$ to $T=10$, the optimization problem is solved with initial states set at $l=1$. The results obtained for the next time step $l=2$ are extracted and fixed at this solution and the prediction horizon is rolled over for the range $T=2$ to $T=11$ and the problem is solved again, this time with the solution for $l=2$ as the initial conditions. This process is repeated until the end of the scheduling horizon $L$. Thus, for a scheduling horizon $L$, the MPC algorithm is given in Algorithm~(\ref{algorithm3}).\\

\begin{algorithm}[H]
\label{algorithm3}
\caption{MPC algorithm for SAEV operation model.}
\KwData{$D_{ij0}$, $U_{i0}^k$, $A_{i0}^{k\theta_i}$ and $\Gamma_0^k$}
\KwResult{the state variables $D_{ijl}$, $U_{il}^k$, $A_{il}^{k\theta_i}$ and $\Gamma_l^k$ for $\forall l \in [0,L]$ }
$l=0$\;
\While{$l \leq L-1$}{
	Set the initial values $d_{ij0}=D_{ijl}$; $u_{i0}^k=U_{il}^k$, $ a_{i0}^{k\theta_i} = A_{il}^{k\theta_i}$  and $\gamma_0^k=\Gamma_l^k$ for the MILP model in the previous section and solve the model with the optimization length $T$\;
	Set the states for next time step: $D_{ijl+1}=d_{ij1}$, $U_{il+1}^k=u_{i1}^k$, $A_{il+1}^{k\theta_i} = a_{i1}^{k\theta_i}$ and $\Gamma_{l+1}^k=\gamma_1^k$\;
  	l=l+1\;
}
\end{algorithm}

\subsection{V2B cost model}
To be able to evaluate the economic value of providing emergency power service via V2B compared to traditional backup generators, the V2B cost $\mathcal{C}_{v2b}$ is expressed as follows:
\begin{equation}
\mathcal{C}_{v2b}= \mathcal{C}_i + f_{out} (\mathcal{C}_e+\mathcal{C}_r)
\label{eq19}
\end{equation}

\noindent where $\mathcal{C}_i$ refers to the total annual cost of installing and operating the V2B system, $\mathcal{C}_r$, refers to the relocation cost per power outage, $\mathcal{C}_e$ refers to the energy cost per power outage and $f_{out}$ is the annual frequency of power outage.

\noindent These costs are calculated as the equations below:
\begin{equation}
	\mathcal{C}_{i}=K\mathcal{C}_v
\label{eq20}
\end{equation}
\begin{equation}
	\mathcal{C}_{e}=(\sigma_t +\omega)q_{v2b}
\label{eq21}
\end{equation}
\begin{equation}
	\mathcal{C}_{r}= T_{relo}\vartheta_c \mathcal{B} \sigma_t
\label{eq22}
\end{equation}

where $K$ refers to the number of vehicles in the fleet, $\mathcal{C}_v$ refers to the annual cost of installing and operating the V2B system per vehicle, $\sigma_t$ refers to the electricity price, $\omega$ refers to the cost of cycling the battery incurred when supplying power to the grid, $\mathcal{B}$ refers to the battery capacity of the vehicle, $q_{v2b}$ refers to the total energy discharged by the SAEV fleet to the building ($= \mathcal{B} * \sum_{k,t} g_t^k$), $ T_{relo}$ is the difference of relocation time between the scenario in which an emergency situation occurs and the one where it does not occur and $\vartheta_c$ is the energy consumed per time step in term of SOC.

\subsection{Model assumptions and limitation}
In our work, the validation of the model is of utmost importance in establishing its credibility and reliability. Despite the inability to validate the model through real-world demonstrations or experiments, its internal consistency and performance in comparison to theoretical expectations can be evaluated. To further verify the influence of model parameters, a sensitivity analysis is performed.

The limitations of the model arise primarily from the assumptions made to simplify the simulations given limited data. One such assumption pertains to the transport simulation, where the city's road network is excluded to reduce simulation time and due to the challenges in accurately representing traffic congestion. Another assumption involves the static transport patterns, which may not accurately reflect the dynamic changes in transport patterns resulting from attracting people who currently use alternative modes.

In terms of energy aspects, the model assumes that the detailed electric grid is not considered and that there are no transmission capacity constraints and transmission losses. This assumption is supported by the expectation that charging stations will be centralized and vehicles will be connected directly to the medium or high-voltage transmission grid rather than the distribution network. Additionally, the detailed charge and discharge behavior of batteries is not considered, and constant charge power and efficiency are assumed at any state of charge. This assumption is deemed acceptable for the purpose of this analysis and is further justified by the simulations' conservative charging power levels.

The design of these charging stations (including location, capacity, etc.) was not considered in this study but is planned as a focus for future work. The numerical assumptions in the case study were supported by references whenever possible. Despite the uncertainty in the results due to the lack of real-world examples of commercial shared autonomous vehicle fleets and the experimental nature of the technology, the internal and theoretical consistency should provide ample validation of the model's analytical capacity for the case study presented.

\section{Case study}
In this section, we present a case study based on realistic travel demand data of the Ile-de-France region in France found in \citep{horl2021open}. A relatively large simulation area consisting of 25 departments of the Ile-de-France region is selected in our case to provide insights on mid to long-distance movements and arrangements for the operation of the SAEV fleet. Smaller or larger geographical regions can be considered with no changes to the methodology proposed. 

In the proposed case study, the duration of each time step taken for the optimization $(\tau)$ is calibrated to represent 6 minutes and we consider 10-time steps as the time horizon for each MPC iteration ($T=10$). This represents a 1-hour optimization horizon for each MILP step. A one-day (24 hours) case is, then, optimized within the MPC framework. This represents the full horizon in our case study. The details of the data used are presented in the following sections covering the transport data and the electric power system data. Each optimization run is solved as a MILP in Julia modeling language by Gurobi solver with version 9.1.2 on a computer with CPU i7-9750H, 8GB RAM, and 6 cores. Based on the obtained results, the mean execution time for a single iteration is approximately 2 minutes and 30 seconds. As such, the overall computational time required for a single day is, on average, around 60 minutes. \footnote{Depending on the instance and the numerical performance of the optimizer used within a single iteration, the computational time may be longer. In our experiments, this computational time for a single day (24 iterations) with the optimality gap for each iteration set to $1 \cdot 10^{-4}$ did not exceed 8 hours. We find this to be acceptable for real-time control since the control action that needs to be implemented is related only to the next time step of one hour (which is calculated within a few minutes) and not the full policy obtained for the whole time horizon.}
\subsection{Transport data}
\begin{figure}[H]
    \centering
    \includegraphics[width=0.8\textwidth]{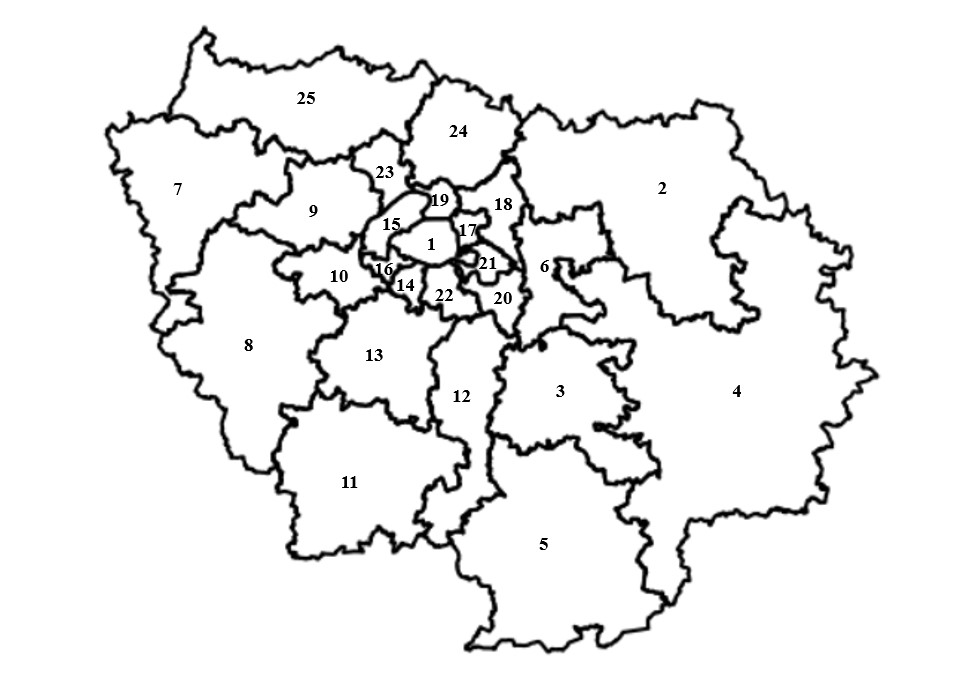}
    \caption{Ile-de-France region and its 25 departments.}
    \label{fig:pic1}
\end{figure}

\begin{figure}[h!]
\centering
    \subfloat[8:00-10:00]{
		\includegraphics[width=0.45\textwidth]{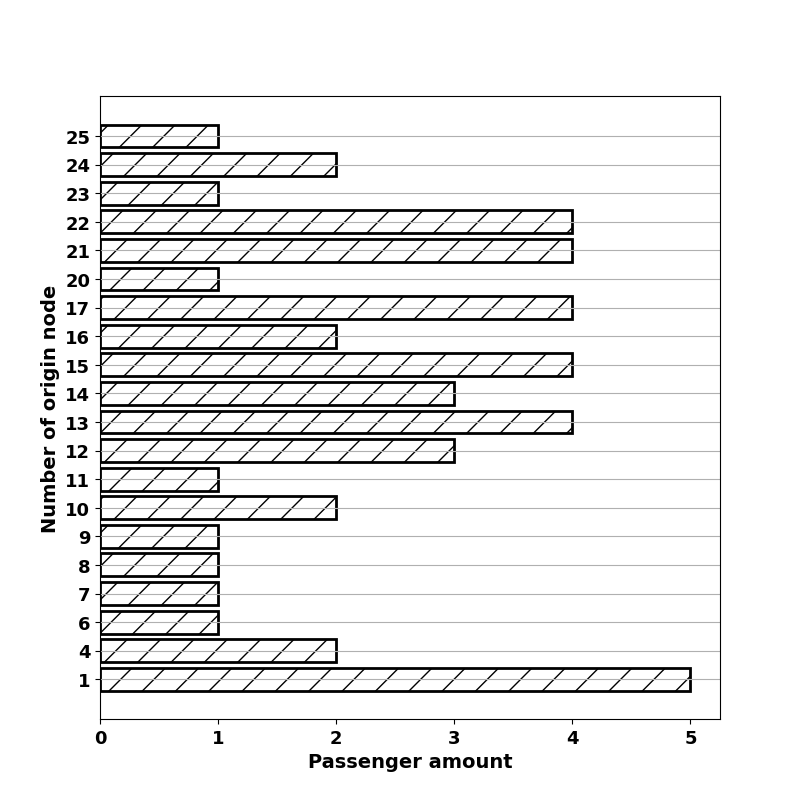}\hspace{-1pt}
		\label{sub51}
		}
    \subfloat[12:00-14:00]{
		\includegraphics[width=0.45\textwidth]{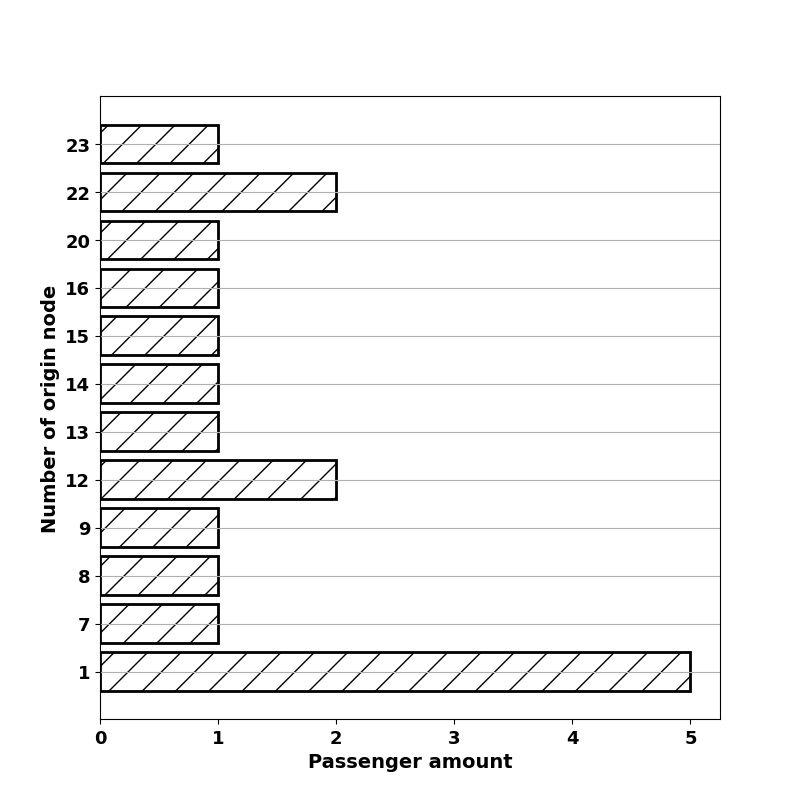}\vspace{2pt}

		}

\vfil
    \subfloat[16:00-18:00]{
		\includegraphics[width=0.45\textwidth]{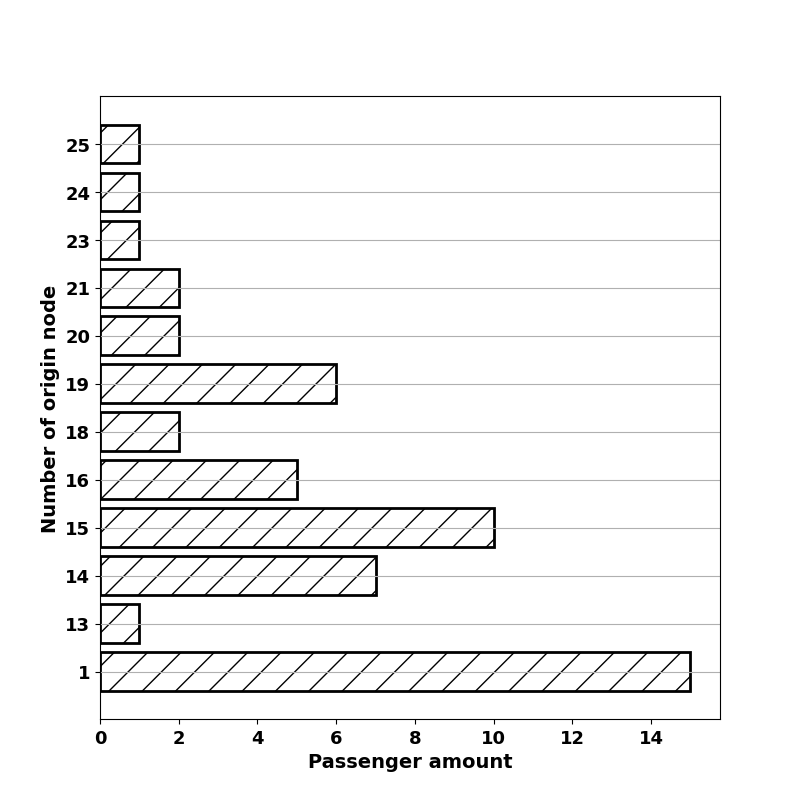}\hspace{-1pt}

		}
    \subfloat[20:00-22:00]{
		\includegraphics[width=0.45\textwidth]{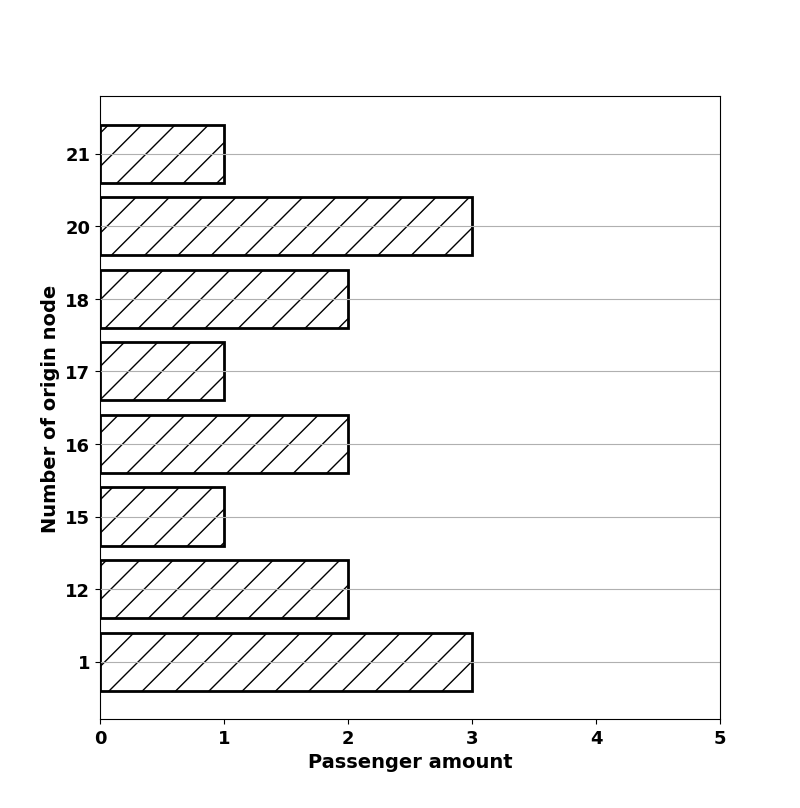}\vspace{2pt}
		\label{sub54}
		}

    \caption{Trip origin distribution in different periods}
    \label{fig:pic2}
\end{figure}

The geographical coordinates of the 25 departments that together comprise the Ile-de-France region are obtained from \citep{Geojson}. Each department is considered a node in our model, either an origin node or a destination node for the trip requests. The map for the Ile-de-France region showing the different departments is shown in Figure  \ref{fig:pic1}, where the numbers refer to department/node numbers.
The trip data in \citep{horl2021open} possesses around 1,040,000 trips in the Ile-de-France area including the origin and the destination coordinates according to the coordinate reference system EPSG:2154, the departure and the arrival time, the travel mode, the travel purpose, among others. Since we consider SAEV services, first, the trips by walk are removed and we only consider trips conducted by both private cars and public transportation. Also, based on our model the trips within one node are not considered. Afterward, the data is separated by different origin-destination (O-D) pairs between the different nodes. The distance between the centroids of the nodes is calculated and, in this work, they are considered as the distances between nodes. Furthermore, these data are divided by departure times at a half-hour frequency in one day. Based on the data, a probabilistic model of passenger arrival times is built. Poisson distribution is chosen to express the number of passengers arriving at each time step \citep{iacobucci2018modeling}. The random distribution of each time step between each O-D pair is obtained by setting the Poisson parameter $\lambda$ equal to the number of trips divided by the number of time steps in each 30-minute period. Finally, a trip distribution is built and around 292 trips are generated for one day as the trips that should be fulfilled with the particular SAEV fleet considered. The distribution of the trips generated along the day is presented in Figure \ref{fig:pic2}.

\subsection{General model parameters}


\begin{table}[H]
\centering
 \begin{tabular}{ccccccccc} 
 Parameter  & $\gamma^{max}$ & $\gamma^{min}$ &$\gamma^{init}$ & $T (hour)$ & $\vartheta_{d}$  & $\vartheta_{c}$ & $\vartheta_{v2b}$ \\\hline
  value& 1	&0.2&	0.8	&1&	0.0092&	0.01&	0.01\\
   & & & & & & & & \\
 Parameter   &$\rho_{1}$ & $\rho_{2}$ & $\eta$ & $K$  & $\omega $ (euro/kWh)&$L$ & $\tau$ (min) \\\hline
  value &	0.01 & 0.001&	0.9&	30	&0.07974&	240&	6 \\
\end{tabular}
 \caption{Values of the parameters used in the case study.}
 \label{table:2}
\end{table}

Table~(\ref{table:2}) summarizes the values chosen for the model parameters. The values for the SOC ($\gamma^{max}$, $\gamma^{min}$ and $\gamma^{init}$) and discharging parameters (discharging efficiency $\eta$ and battery cycling cost $\omega$) are comparable to the values used in references \citep{iacobucci2019cascaded} and \citep{iacobucci2019optimization}. For the vehicle characteristics, parameters roughly similar to the catalog values of the Tesla Model S are taken in our model with 85kWh battery capacity, 10 hours for the slow charge, and a cruising range of 652 km. The electricity price is taken as 0.1292 euro/kWh, which is the average household electricity cost for the Ile-de-France region in 2020 between July and December, as found in \citep{Electricityprices}.

\subsection{Parameters for the emergency case}
We consider the energy demand during a power outage of a large hospital in the Ile-de-France region as an emergency case. The work in \citep{bawaneh2019energy}  has shown the annual energy intensity of healthcare facilities in France to be approximately $228.2 kWh/m^2$. We consider a case study of a hospital in the west of Paris with an area of approximately $120,000 m^2$. By averaging the annual energy cost, the energy demand per time step $Q_d$ can be obtained as $312.6 kWh/time step$. This translates to $Q_d= 3.678$ as a function of SOC to keep the consistency of units in the model. We assume a single SAEV fleet operator with a fleet size of 30 vehicles and a battery charging rate $\vartheta_c=0.01$. Since a power outage for a hospital of this size would not be possible to satisfy solely with one SAEV fleet with the size proposed, we assume with no lack of generality of the proposed approach that the hospital is connected to six 500 kW backup generators and we focus on comparing the cost of adding one more generator to the cost of utilizing V2B power provided by the SAEV fleet to meet the remaining 126 kW demand. The work in \citep{grosskopf2006assessing} has presented that in 5 years term, a 500 kW generator has a planned value (PV) cost of \$74,264 (approximately 13,367 euros/year), which will be used for cost comparisons and analysis.

Finally, referring to \citep{sarparandeh2017pricing} and the assumptions presented above, the values of parameters related to V2B costs are as shown in Table \ref{table:3}.
\begin{table}[H]
\centering
 \begin{tabular}{c|c|c|c|c|c|c} 
 parameter & $K$&$\mathcal{C}_v$&	$\sigma_t$&	$\omega$&	$\mathcal{B}$	&$\vartheta_c$ \\ \hline
 value&30&	45 euro/vehicle/year &	0.1292 euro/kWh	&0.0797 euro/kWh	&85 kWh	&0.01\\

\end{tabular}
 \caption{Values of parameters related to cost calculation.}
 \label{table:3}
\end{table}
\section{Results}

In this section, we present the optimization results based on the case study proposed. First, in Section \ref{normal_results_section} the results under normal conditions (no power outages) are discussed in detail. In Section \ref{emergency_case_section} the results of the emergency case are presented, including the impact of the outage on the SAEV fleet performance, the capacity to respond to the power demand, and the cost evaluation. Finally, a detailed sensitivity analysis is conducted to evaluate the impact of considering different outage locations and outage durations on the performance of the SAEV fleet and on the resilience response for the critical building.

\subsection{Normal operation conditions}
\label{normal_results_section}
First, the results of the optimization problem under normal operating conditions are presented without considering power outages. A system with 30 vehicles representing the SAEV fleet and 25 nodes is considered for a 1-day horizon. We evaluate a number of performance measures related to the case study under normal conditions to establish the capabilities of the model and to provide a baseline for comparison under different scenarios. 


\begin{figure}[h!]
    \centering
    \subfloat[Impact of varying the number of passengers in the system on total passenger waiting time.]{
		\includegraphics[width=0.6\textwidth]{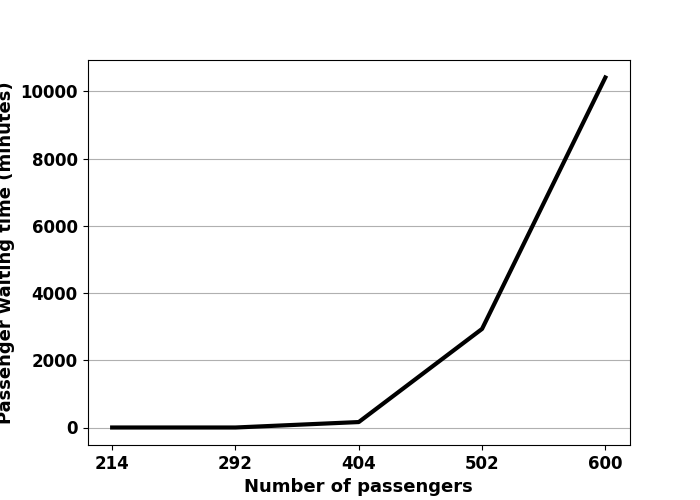}\hspace{5pt}\vspace{1pt}
		\label{fig:pic3a}}

	\subfloat[Impact of varying the number of passengers trip requests on total vehicle relocation times.]{
	\includegraphics[width=0.6\textwidth]{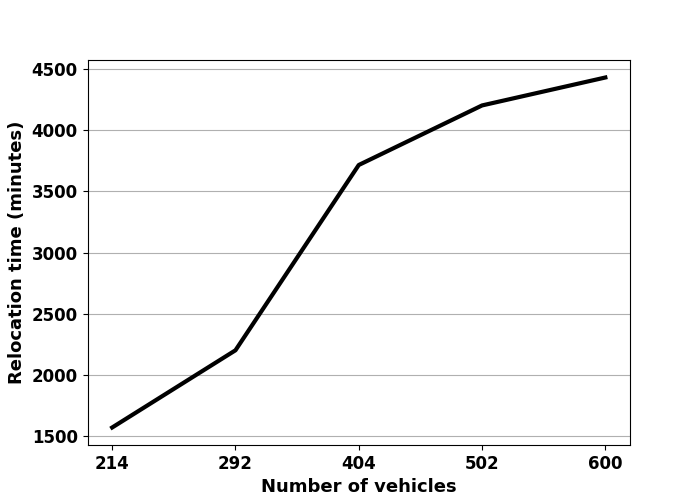}\hspace{5pt}\vspace{1pt}
	\label{fig:pic3b}}
    \caption{Impact of varying the number of passengers on the performance of the SAEV fleet.}

\end{figure}

Figure (\ref{fig:pic3a}) shows the impact of varying the number of passengers requesting mobility services on the total passengers waiting time. As shown in the figure, 30 vehicles can satisfy the transportation requirements of 300 passengers with no passengers waiting. Therefore, a lower number of passengers trip requests in the network can still be satisfied with no waiting time. However, adding 100 new passengers leads to a rapid increase in the passengers waiting time. The total waiting time in the system increases from 0 to 162 minutes for around 400 passengers, then to 2934 minutes for 500 passengers, and then to 10,406 minutes for 600 passengers, as a result of the accumulation of waiting passengers. 

The impact of increasing the number of passengers in the system on the total relocation time of the vehicles is shown in Figure~(\ref{fig:pic3b}). It is seen that, for optimal results, the case with 300 passengers requires 2202 minutes for vehicles to relocate. However, with 100 passengers less in the system, the relocation time drops to about 1572 minutes. On the other hand, increasing the number of passengers in the system increases the total relocation time. This increase is approximately logarithmic in nature (the rate of increase is declining as the number of passengers increases). For 400 passengers the total relocation time of the vehicles is around 3714 minutes (a 1512 minute increase or 69 \% increase compared to the case with 300 passengers). Further increasing the number of passengers to 500 passengers leads to a total relocation time of 4200 minutes (a further 13\% increase compared to 400 passengers). Finally, for a system with 600 passengers, the total relocation time reaches 4428 minutes.

These results can be explained with respect to the passengers’ total waiting times. As the objective of the model is to minimize the passengers waiting time, for the scenarios of 300 passengers or less, the fleet is capable of satisfying all travel demands and there are no passengers waiting times. In these cases, the reduction in the number of passengers directly leads to a reduction in the total relocation time of the vehicles, since the reduced number of passengers reduces the overall mobility needs in the system. Reciprocally, increasing the number of passengers above 400 requires more relocations to reach all passengers. The relocation of the SAEV vehicles increases as the number of passengers increases, however, since the number of vehicles remains constant, the increase in the relocation time approaches a maximum, after which, regardless of the number of passengers in the system, no further relocation can be made in order to satisfy the mobility demands.

\begin{figure}[h!]
    \centering
    \subfloat[ Impact of varying the number of vehicles of the SAEV fleet on the total passenger waiting time.]{
		\includegraphics[width=0.6\textwidth]{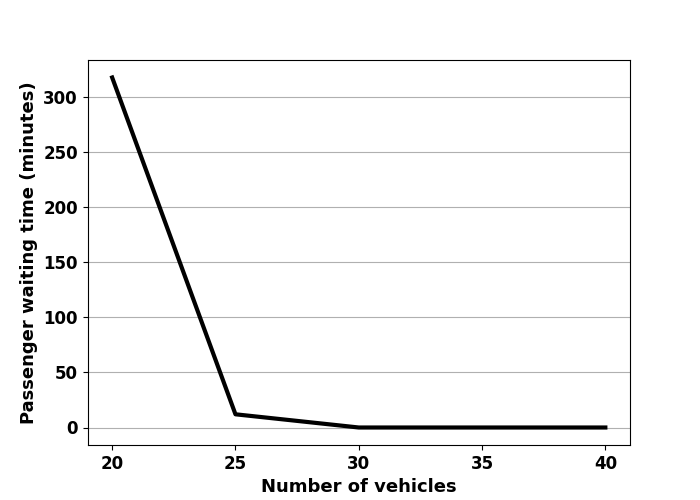}\hspace{5pt}\vspace{1pt}
		\label{fig:pic4a}}

	\subfloat[Impact of varying the number of vehicles of the SAEV fleet on the total relocation time.]{
	\includegraphics[width=0.6\textwidth]{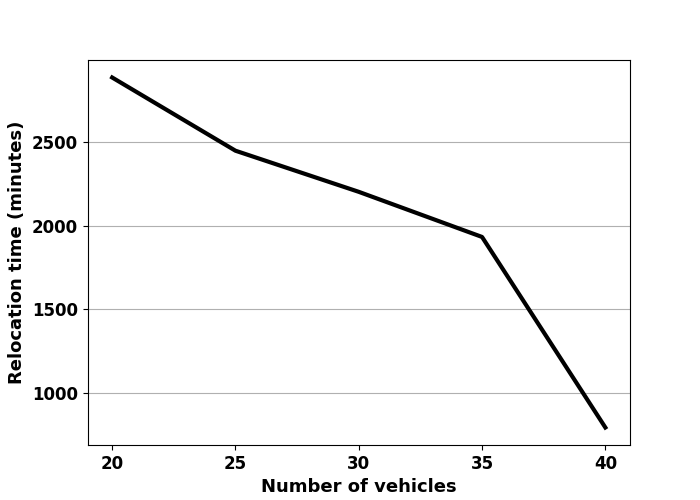}\hspace{5pt}\vspace{1pt}
	\label{fig:pic4b}}
    \caption{Impact of varying the number of vehicles of the SAEV fleet on the performance of the fleet.}

\end{figure}

Next, we study the impact of varying the number of vehicles in the SAEV fleet on the performance of the system. The number of passengers in the system is fixed at 292, and the impact of changing the number of vehicles on the passengers total waiting time and on the relocation time is shown in Figure~(\ref{fig:pic4a}) and Figure~(\ref{fig:pic4b}), respectively.

Based on the same case with 30 vehicles and 292 passengers, Figure~(\ref{fig:pic4a}) shows that increasing the number of vehicles beyond 30 does not appear to have any impact on the passengers total waiting time as the optimal solution results in a passenger waiting time of 0 in all these cases. In contrast, reducing the fleet to 25 vehicles has a slight impact on increasing passengers waiting time to 12 minutes. Further reducing the fleet size to 20 vehicles considerably increases the passengers waiting time to 318 minutes. This rapid increase results from the loss of the transportation capacity of the fleet and would increase even further for smaller fleet sizes. 

The results in Figure~(\ref{fig:pic4b}) illustrate the impact of varying the number of vehicles on the total relocation time in the network. It is seen that, compared to the case of a fleet size of 30 vehicles, a reduction of the fleet size to 25 vehicles increases the relocation time from 2202 to 2248 minutes (an 11\% increase compared to the case with 30 vehicles). A fleet size of 20 vehicles would further increase the relocation time to 2886 minutes (a 19\% increase compared to the case with 25 vehicles). It is clear that reducing the fleet by 10 vehicles, or in other words reducing the transportation capacity of the fleet, will create more relocation orders. Conversely, as expected, increasing the fleet size reduces the total relocation time in the system. For the case of 35 vehicles, the total relocation time is reduced from 2202 minutes to 1932 minutes (a 19\% decrease compared to the case with 30 vehicles). Furthermore, a fleet size of 40 vehicles reduces the relocation time to 792 minutes(a 59\% decrease compared to the case with 35 vehicles). In this case, the significant decrease in relocation time is a result of more options for dispatching vehicles provided by the greater fleet size.

\subsubsection*{Evaluating the charging behavior of SAEV fleet under normal operating conditions:}

\begin{figure}[h!]
\centering
    \subfloat[8:00-11:00]{
		\includegraphics[width=0.45\textwidth]{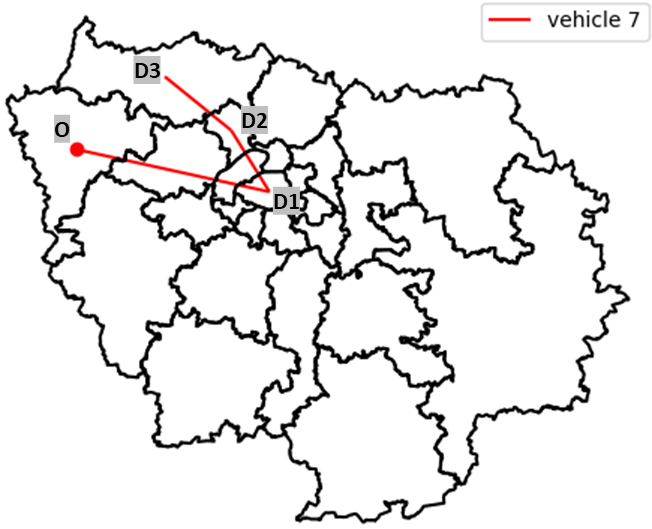}\hspace{-1pt}
		\label{sub51}
		}
    \subfloat[11:00-14:00]{
		\includegraphics[width=0.45\textwidth]{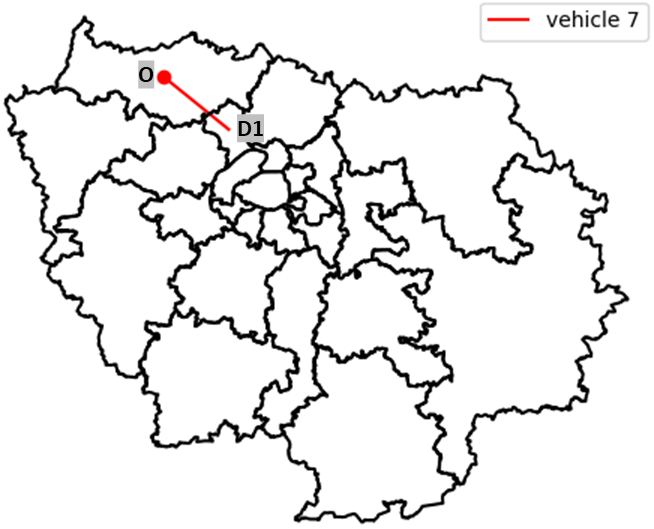}\vspace{2pt}

		}

\vfil
    \subfloat[14:00-17:00]{
		\includegraphics[width=0.45\textwidth]{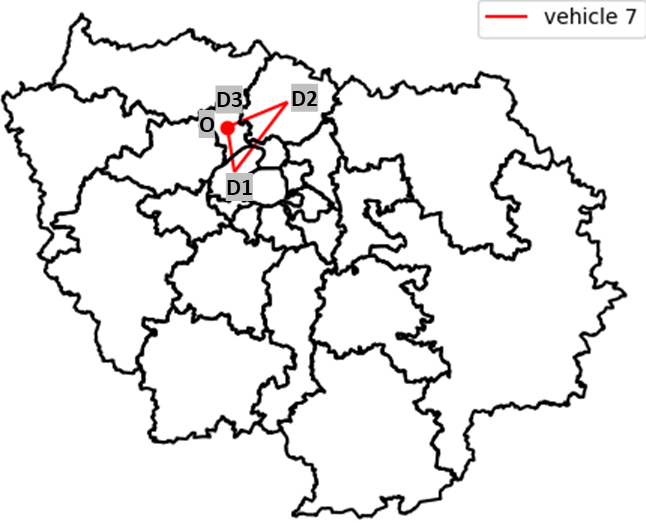}\hspace{-1pt}

		}
    \subfloat[17:00-20:00]{
		\includegraphics[width=0.45\textwidth]{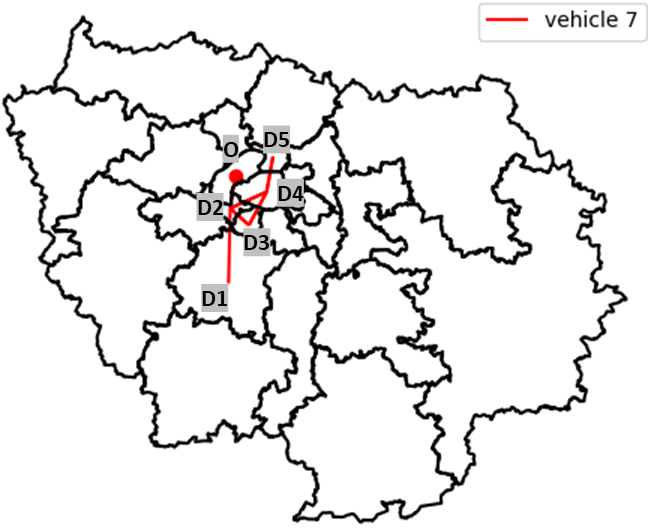}\vspace{2pt}
		\label{sub54}
		}

    \caption{Trajectory of vehicle 7 throughout the day under normal operational conditions.}
    \label{fig:pic5}
\end{figure}

Since SAEV fleet performance is dependent on the charging state of the vehicles, it is important to analyze the results regarding the energy consumption and recharging behavior of the vehicles. We take the case of 30 vehicles and 292 passengers in one day  as a base scenario. First, an example of vehicle trajectories throughout a day is shown in Figure~(\ref{fig:pic5}).

Figure~(\ref{fig:pic5}) illustrates  the trajectory of vehicle 7 as an example. The figure breaks down the trajectory of the vehicle during 3 hours for the 12-hour duration between 8h00am and 8h00pm (Figures~\ref{sub51} to \ref{sub54}) to illustrate the dynamic optimization results of the model. Vehicle 7 undertakes transport demand from 8:00-20:00 and the nodes in the figures represent the departure node at the beginning of every 3 hours time period. 

Generally, the trips of most of the vehicles of the fleet are within the central area of the Ile-de-France region and vehicles might pass by the same region several times during the day either to relocate or to transport passengers. A detailed illustration of the trajectories taken by all the vehicles of the SAEV fleet throughout the day as obtained by the optimal solution can be seen in Figure~(\ref{appen:pic1}) in \ref{appen1}.

\begin{figure}[h!]
    \centering
    \includegraphics[width=1\textwidth]{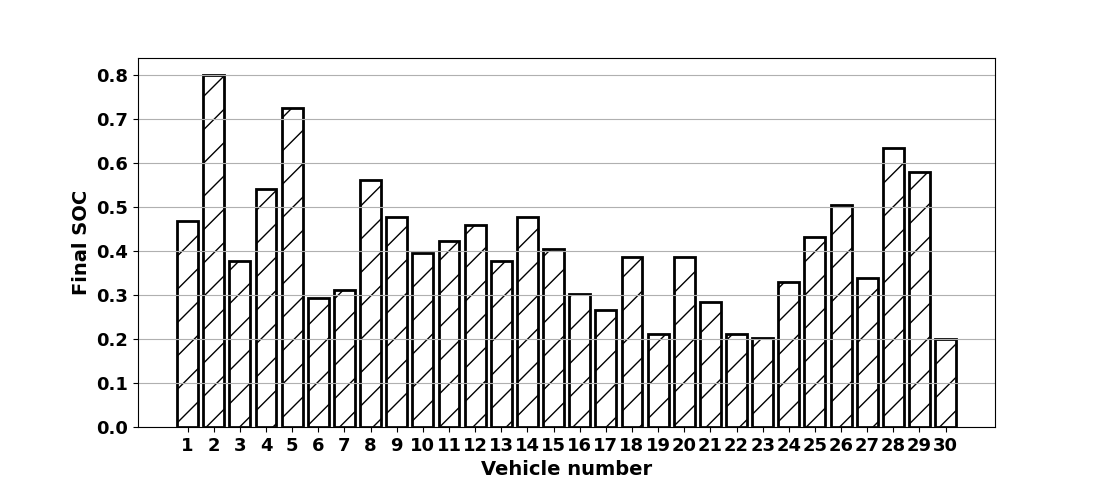}
    \caption{Final SOC of the vehicles at the end of the one-day time horizon. Excess energy above 0.2 represents potential energy that can be used to provide emergency V2B services.}
    \label{fig:pic6}
\end{figure}

 In this case study under normal operating conditions, most of the vehicles do not need to be charged. In fact, it can be noted that with the initial SOC 0.8 and in a one-day time horizon, the energy requirement of the system is low, i.e., vehicles are sufficiently charged to satisfy the trip requests without significant charging needs (with the exception of vehicle number 30 which charges for around 30 minutes at node 22 from 21h36 to 22h06). The final SOC of all the vehicles in this case is presented in Figure~(\ref{fig:pic6}). It is shown that, at the end of the day, the SOC of the vehicles varies between 0.2 (vehicle 30) to 0.8 (vehicle 2), with none of the vehicles below the minimum SOC permitted. Any vehicle with a SOC higher than 0.2, therefore, has the potential to provide V2B services if needed with no impact on its charging needs within a day.


\subsection{Power outage emergency case}
\label{emergency_case_section}

In the emergency case, it is first assumed that the power outage occurs from 17:30 to 18:30
at node 1. As before, the fleet size is 30 vehicles and the passenger demand is set to 292. The case study is solved with the proposed model, considering the need for vehicles to provide electricity to critical buildings when an outage occurs. 
\subsubsection{Evaluating the performance of the SAEV fleet and critical building under emergency power outage:}

\begin{figure}[h!]
    \centering
        \subfloat[]{
		\includegraphics[width=0.5\textwidth]{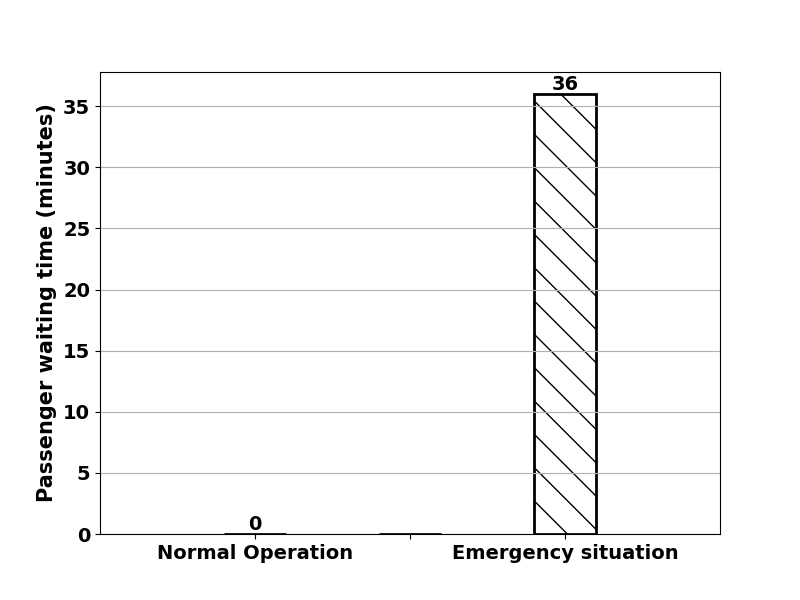}
		}
    \subfloat[]{
		\includegraphics[width=0.5\textwidth]{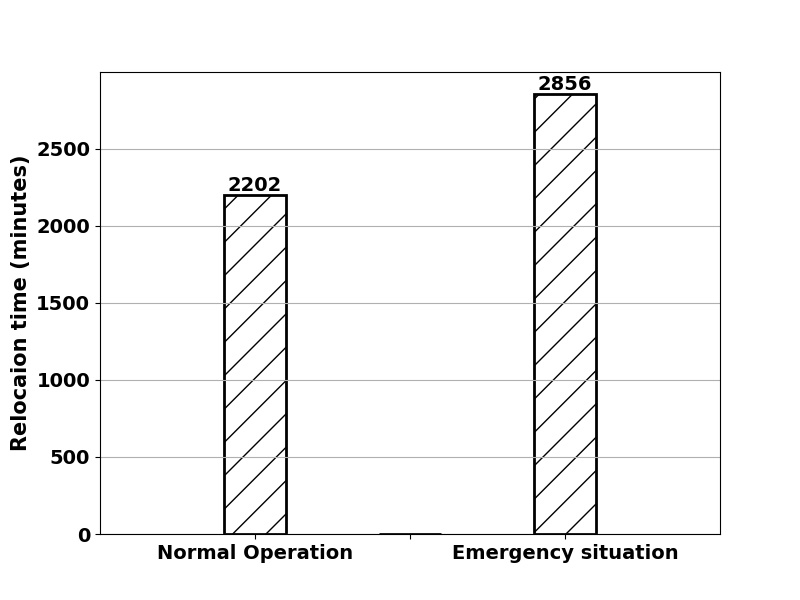}
		}
    \caption{Comparison of emergency situation and normal situation}
    \label{fig:pic7}
\end{figure}

First, a comparison of the total passenger waiting time and relocation time in the system between the normal operating scenario and the emergency case is presented in Figure~(\ref{fig:pic7}). As seen in the figure, both passenger total waiting time and relocating time increase in the power outage case due to the V2B contribution of the SAEV fleet towards the power needs of the hospital. This results in fewer vehicles available to carry passengers and in higher relocation to arrive at the hospital location.

The passengers total waiting time increased in this case to around 36 minutes. However, it would be important for the fleet operator to evaluate this loss in transportation performance with respect to the added resilience towards the power outage in the critical building. For this, we perform a cost analysis to evaluate the cost of providing the needed power via the SAEV fleet, compared to installing additional backup generators for the critical building. 

From the optimization results, we can calculate the difference of the vehicles relocation time under normal operational conditions and in the emergency situation $ T_{relo} =654 \text{min}$, as well as the total energy discharged by V2B $q_{v2b} =139.74$ kWh. Then, based on Eqs.~(\ref{eq19}-\ref{eq22}) and the data in Table~(\ref{table:3}), we can obtain the  cost of installing and operating a V2B system for the critical building  $\mathcal{C}_i=1350 \text{\euro}/year$, the energy cost $ \mathcal{C}_e=29.19 \text{\euro}$ and the relocation cost $\mathcal{C}_r=71.82 \text{\euro}$. Thus, by comparing the final annual cost, we can calculate that introducing another 500 kW generator is better than introducing a V2B system only when the outage frequency $f_{out}\geq 121$. That is to say only if the power outage in this case occurs more than 121 times per year.

Under the emergency case, the discharging behavior appears as a result of V2B demand. The discharging behavior of the SAEV fleet is described in Table~(\ref{table:5}) in terms of the total discharge amount of each vehicle during the outage duration. Only the vehicles that have contributed to V2B services during the power outage are reported in the table.
 
\begin{table}
\centering
\resizebox{\textwidth}{!}{
 \begin{tabular}{p{2.5cm}|c|c|c|c|c|c|c|c|c|c|c } 
 Vehicle number & 1&2&3&6&10&11&12&13&14&16&17\\\hline
 Total discharge amount/SOC & 
 0.05 &0.1&0.1&0.1& 0.07 &0.02&0.05&0.0676&0.02&0.09&0.1\\\hline
 Vehicle number & 18&19&21&22&23&24&25&26&28&29&30\\\hline
  Total discharge amount/SOC & 
0.0388 &0.1&0.0493&0.03& 0.07 &0.02&0.00944&0.1&0.00944&0.09&0.0796\\\hline
\end{tabular}}
 \caption{Vehicles total discharging amount to the critical building during the outage duration.}
 \label{table:5}
\end{table}

The results summarized in Table~(\ref{table:5}) show that 23 vehicles among 30 vehicles in the fleet have discharged during the emergency. Among all the vehicles that have contributed to the V2B service, 7 vehicles (vehicles 2,3,6,17,19,26, and 27) continuously discharge during the whole outage duration with the maximum rate (0.01), whereas the other vehicles provide a lower discharging rate at different times during the outage duration. This is because the V2B requirement of the critical building does not require maximum discharging energy from all the vehicles at the same time and, instead, for the optimal solution, some vehicles would not connect to the building in some periods and would instead provide mobility services.

\begin{figure}[h!]
    \centering
        \subfloat[Comparison of trajectory of vehicle 7.]{
		\includegraphics[width=0.5\textwidth]{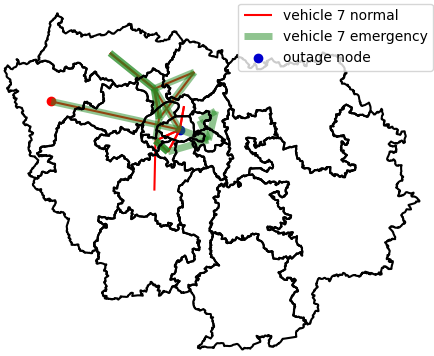}\hspace{-1pt}
		\label{fig:pic8_1}
		}
    \subfloat[Comparison of trajectory of vehicle 8.]{
		\includegraphics[width=0.5\textwidth]{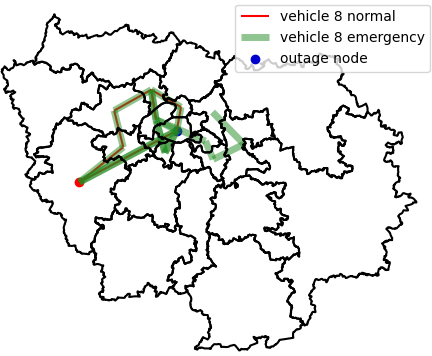}\hspace{-3pt}
		\label{fig:pic8_2}
		}

    \subfloat[Comparison of trajectory of vehicle 17.]{
		\includegraphics[width=0.5\textwidth]{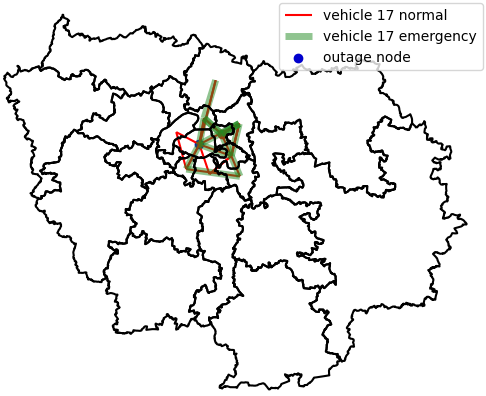}
		\label{fig:pic8_3}
		}
    \caption{Comparison of trajectories under normal conditions and during the emergency power outage for selected vehicles.}

\end{figure}

Next, we investigate the variations in the vehicles' trajectories due to the presence of the emergency case. Several typical variations are shown in Figure~(\ref{fig:pic8_1}), Figure~(\ref{fig:pic8_2}) and Figure~(\ref{fig:pic8_3}). As shown in Table~(\ref{table:5}), several vehicles are required to provide V2B power to the critical building and, therefore, this leads to a change in trajectory for most of the vehicles compared to the normal case. This is because the optimization problem seeks to optimize transportation needs as some vehicles are taken out of the mobility system to provide the required power to the building.

Several typical changes in vehicle trajectories are shown in  Figure~(\ref{fig:pic8_1}), Figure~(\ref{fig:pic8_2}) and Figure~(\ref{fig:pic8_3}). For example, Figure~(\ref{fig:pic8_1}) illustrates the case where significant changes in the vehicle trajectory occur during the day. Some trips are abandoned while some other trips are undertaken as a result of the fleet behavior change. This is one of the most common cases: the global fleet movement is changed by the power emergency request, therefore, the optimal solution leads to a remarkable change in trips undertaken between vehicles as is the case for vehicle 7. Another possible behavior is shown in  Figure~(\ref{fig:pic8_2}), where vehicle 8 may follow an overall similar trajectory, but with a significant extension as a result of the need to relocate to the power outage node. Note that vehicle 8 does not pass by node 1 in normal operation, but relocates to node 1 during an emergency. Finally, Figure~(\ref{fig:pic8_3}) shows a relatively rare situation: Vehicle 17 reduces its trajectory compared to the no emergency case. This can be explained with regards to the fact that vehicle 17 contributes the most to the power need of the critical building as seen in Table~(\ref{table:5}) and, therefore, it is primarily excluded from providing mobility services.

\subsubsection{Sensitivity analysis for power outage time and location}

\begin{figure}[H]
    \centering
        \subfloat[Impact of power outage start time on total passengers waiting time.]{
		\includegraphics[width=0.6\textwidth]{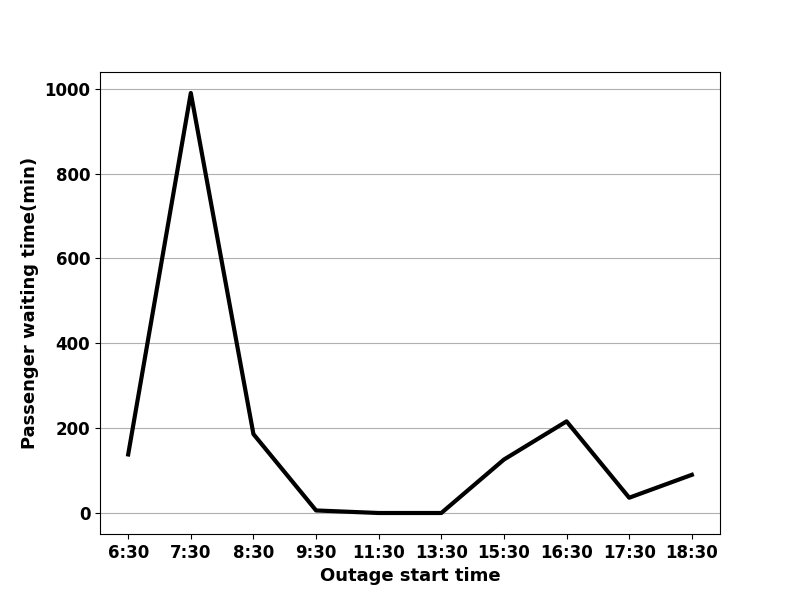}
		\label{fig:pic9_1}
		}

    \subfloat[Impact of power outage start time on total fleet relocation time]{
		\includegraphics[width=0.6\textwidth]{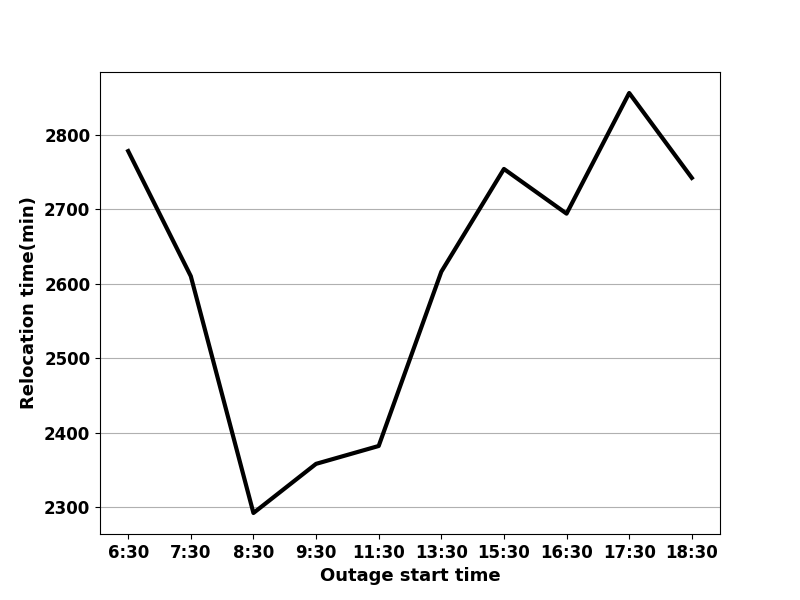}
		\label{fig:pic9_2}
		}

    \caption{Sensitivity results for power outage start time.}
\end{figure}

We further investigate power outage scenarios under different power outage start times and locations. Different power outage scenarios are explored with the same 1-hour outage length, the number of vehicles and the number of passengers demand for one day. The impact of solving the optimization problem while varying the power outage location and start times is studied. Figure~(\ref{fig:pic9_1}) illustrates the sensitivity results for the total passengers waiting time in the system, and Figure~(\ref{fig:pic9_2}) illustrates the results for the impact on the total relocation time of the SAEV fleet. 

As seen in Figure~(\ref{fig:pic9_1}), passengers total waiting time caused by the emergency power demand achieves two peaks when the power outage occurs between 7:30- 8:30 and between 16:30-17:30, which refers to the major peak hours in the morning and evening as shown in Figure~(\ref{fig:pic2}). Meanwhile, compared to 216 minutes of passengers total waiting time caused by an outage at 16:30-17:30, an outage at 7:30-8:30 leads to 990 minutes of total passengers waiting time. This can be explained by the difference in the origin nodes in the two situations. During the morning, most passengers arrive at nodes 10, 12, 13, and 14 which are further from the outage node in this case (node 1). This means that it is more difficult for the SAEV fleet to both respond to the emergency power need and satisfy the transportation demand. This is not the case for the outage during evening time as the majority of the passengers arrival is at node 1 (the location of the power outage) and the surrounding nodes 15, 16, and 19. 

\begin{figure}[h!]
    \centering
    \includegraphics[width=1\textwidth]{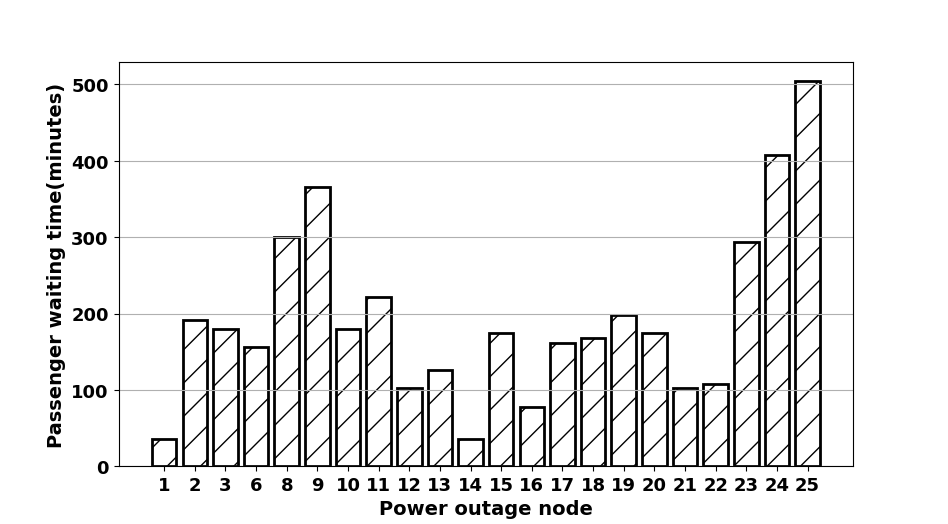}
    \caption{Impact of varying the outage location on the total passengers waiting time in the system.}
    \label{fig:pic10}
\end{figure}

Furthermore, the impact of different power outage locations is also investigated. As before, the results shown are based on the 30 vehicles and 292 passengers case in the one-day operational horizon. The power outage time is set to 17:30-18:30. The total passengers waiting time in the system is considered as a primary indicator to show the impact on the mobility system. The impact of the changes in power outage nodes on the total passenger waiting times is shown in Figure~(\ref{fig:pic10}).

First, it should be noted that if the power outage occurs at node 4,5, or 7, the model will provide an infeasible result as the traveling time from most nodes to these nodes are longer than the prediction horizon. This is a limitation of our model and can be avoided by ensuring that all the trip times between nodes are beneath the prediction horizon. However, this can be considered acceptable for our case study since, as shown in Figure~(\ref{fig:pic2}), very little passenger demand originates at these furthest nodes and, therefore, would not significantly impact the results. Furthermore, the infeasibility of the solution in these cases can be practically understood as the incapacity of the SAEV fleet to respond to these emergency situations given their original locations at the time the outages occurred. This could be solved by adding more vehicles in the fleet, or explicitly placing part of the fleet closer to these remote areas.

For the rest of the outage scenarios, the impact on the total passengers waiting time is driven by a complex interaction between the node location and the number of passengers arriving and departing at this node. Generally, as two major passengers destinations are node 1 and node 14, the passenger waiting time is the lowest for an outage that occurs at these nodes. For further outage locations, the waiting time can become significantly longer, reaching levels higher than 300 minutes and up to 500 minutes for nodes 8, 9, 23, 24, and 25.

\begin{figure}[H]
    \centering
    \includegraphics[width=1\textwidth]{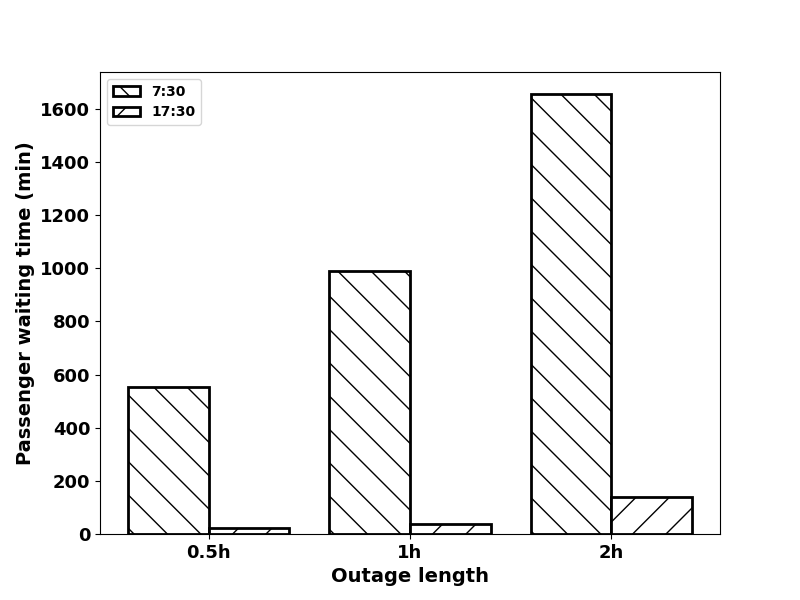}
    \caption{Impact of the outage length.}
    \label{fig:pic11}
\end{figure}

Additionally, an investigation was conducted to examine the impact of outage length on the system. The scenario under consideration assumes that the outage occurs at node 1. To capture the effect of varying outage durations, simulations were carried out with the outage start times set at both 7:30 and 17:30. The simulation results are presented in Figure~(\ref{fig:pic11}).

According to the figure above, the impact of longer power outage durations on passengers total waiting times is not simply additive. The passenger total waiting time resulting from longer power outages is significantly greater than the sum of the waiting times from two independent shorter power outages. This is because more prolonged power outages often result in vehicles within the system being unable to meet passenger demand for a more extended period of time, leading to an accumulation of unmet demand. Additionally, during periods when the system is more resilient against power outages (such as the period between 17:30 and 19:30 in the figure), the impact of longer power outages on the system is not as significant.

\begin{figure}[H]
    \centering
    \includegraphics[width=1\textwidth]{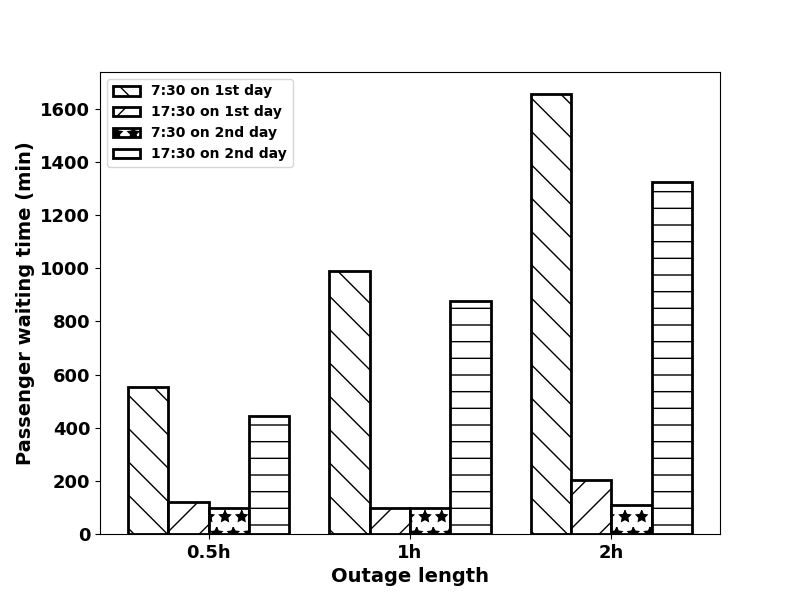}
    \caption{2-day optimization result.}
    \label{fig:pic13}
\end{figure}

Moreover, a two-day case simulation has been developed, utilizing outage node 1, and different outage lengths have been estimated. The simulation parameters dictate that the outage start time be set at either 7:30 or 17:30 on the first or second day. It is also assumed that vehicles will be charged to the initial state after each day. The simulation results have been presented in Figure ~(\ref{fig:pic13}).

Given that we assume all vehicles will be charged to their initial state at the end of each day, the primary impact of the operations during the first day on those of the second day is on the locational distribution of the vehicles. Therefore, we observe that in the case of considering two days, the total passenger waiting time is 96 minutes for the scenario with a power outage at 17:30 on the first day, which is higher than the passenger waiting time of 36 minutes in the one-day case with the same power outage time. The analysis shows that the extra 60 minutes of passenger waiting time all occurred on the second day. This is due to the late timing of the power outage which leads to an abnormal concentration of vehicles  to the power outage node at the end of the first day. Conversely, if the power outage occurs at 7:30 on the first day considering the two-day horizon, the passenger waiting time was the same as the previously found results considering a one-day horizon, and there was no passenger waiting on the second day at all. This is because the subsequent operation on the first day after the outage caused all vehicles to return to a similar distribution among the nodes and therefore not impacting the operations on the second day. 

For cases where a power outage occurs on the second day, the change in customer total waiting time is determined solely by the relationship between the second day's passenger arrival times and the outage time. The impact of different outage times varies due to the difference in passenger distribution between the first and second day, which can be referred to in the explanation of Figure ~(\ref{fig:pic9_1}). 

\begin{figure}[H]
    \centering
    \includegraphics[width=0.7\textwidth]{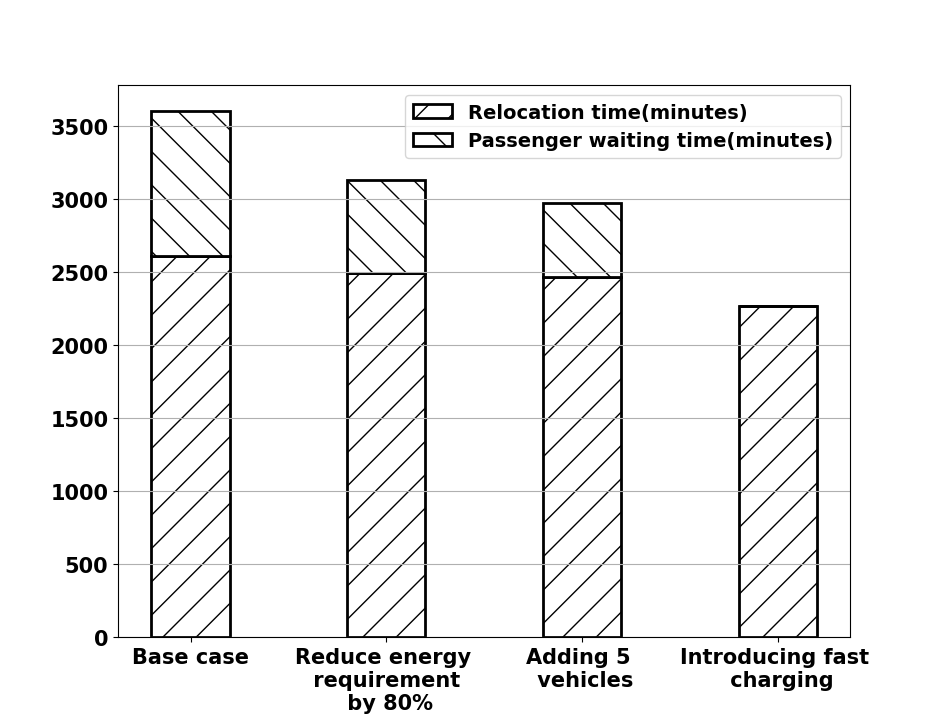}
    \caption{Impact of different operational strategies on the total passengers waiting time and relocation time during a critical building power outage.}
    \label{fig:pic12}
\end{figure}

Finally, noticing that, for the case study presented, if the outage occurs at 7:30-8:30 at node 1, the passenger total waiting time will become significantly high, we explore different strategies to reduce the passenger waiting time. Particularly, we investigate the impact of reducing V2B energy needed, adding vehicles to the fleet, and introducing fast charging. The results of implementing these strategies on the optimization results in terms of total passengers waiting time and relocation time are shown in Figure~(\ref{fig:pic12}).
As shown in the results presented in Figure~(\ref{fig:pic12}), reducing the energy demand in emergency situations by 80\% (for example by relying on additional backup generators) or adding 5 more vehicles to the fleet can provide a reduction in passengers total waiting time. However, the adoption of fast charging ($\vartheta_c=\vartheta_{v2b}=0.1$) reduces passengers total waiting time effectively while avoiding the need to increase the number of vehicles in the fleet or increasing the investment costs in backup generators. The reduction of the total passengers waiting times in this case is 990 minutes compared to the base case, which is a significant improvement compared to the other scenarios. Similarly, the reduction of total relocation times is 340 minutes compared to the base case, which is a slight improvement compared to the other strategies.  Therefore, an installation of fast charging can be regarded as a considerable means of improvement for V2B service resilience in an emergency.
\subsection{Summary of the main results, analyses, and decision-making insights:}

\noindent In the previous sections, we show the capability of the proposed modeling and optimization framework in addressing the problem of coordinated autonomous traffic assignment and V2B response to critical energy emergencies of an SAEV fleet. The results are discussed in detail to show the large number of insights that can be obtained using the proposed dynamic optimization framework to treat this problem. Here, we attempt to summarize the main insights as seen in the results that highlight the usefulness of the proposed approach for practical decision-making.

The main results and insights of this part can be summarized as follows:

\begin{itemize}
  \item The model shows a logical relation between the decrease in the SAEV fleet capacity (due to a higher number of passengers of lower numbers of vehicles) and the performance of the SAEV fleet: the increase of arriving passengers beyond the fleet capacity results in a rapid decrease in the fleet performance, as measured by the total waiting times of the passengers and relocation times of the vehicles. The relocation time follows a logarithmic increase as the capacity of the fleet is reached, i.e. beyond a certain level no further increase in relocation time occurs and the worsening in performance is strictly passed to the passengers waiting times.
    \item Based on the sensitivity analysis of the transportation capacity parameters, 30 vehicles, and 292 passengers case study is chosen as the basis for analyzing the results of the model. Further examination shows that, in this case, the SAEV fleet has enough capacity to respond to the transportation demand. Moreover, an investigation of the energy status of the vehicles shows that only vehicle 30 required charging and most of the vehicles had around 30\% - 50\% available SOC throughout the day, which highlights their potential for providing V2B services.
  \item An outage is assumed to be at node 1 from 17:30-18:30 and the V2B potential is evaluated. As a result of the V2B emergency response, the passenger total waiting time increases by 36 minutes, and the relocation time increases by 654 minutes compared to the base case (no power outage). A further cost evaluation proves that compared with adding another backup generator, a V2B service provides a lower cost per year unless the same outage occurs more than 121 times per year.
  \item  A sensitivity on the outage start time is evaluated and two peaks of passenger total waiting time are found (990 minutes waiting time if the outage occurs at 7:30-8:30 and 216 minutes waiting time if the outage occurs at 16:30-17:30). These two peaks result from the higher transport demand in these time periods (peak hours) which highlight the need to pay special attention to reinforce the service of the system during these periods or to maintain back-up generators ready  around those times if vehicles are not available to respond to a power emergency. 
  \item The impact of the outage node on the performance of the SAEV fleet is also investigated and the passengers destinations appear to be a major factor in impacting the performance of the SAEV fleet during a power emergency. A careful study of passengers traveling patterns can, therefore, better prepare the policy-makers for assessing the risk of not having enough vehicles to respond to a critical energy emergency in further locations.
  \item Finally, several strategies to reduce the impact of SAEV V2B response to a power outage on the passengers total waiting times are investigated. The results show that adopting fast charging at the critical building's location significantly reduces the passengers total waiting time as fewer vehicles are needed to provide power at a given instant. 
\end{itemize}

It should be noted that the analysis provided in our work focuses on scenarios of a single power outage at a time, considering different outage lengths and locations. In these cases, it is shown that SAEVs provide a viable and more economical alternative to respond to a critical building's energy needs, compared to traditional backup generators. The main hypothesis for this choice is that full power outage situations are, often, rare events (particularly in developed economies where the security of electricity supply is of high importance to infrastructure management). A direct extension of this analysis could consider the case of several outages occurring simultaneously on different nodes, or successively throughout the day to evaluate the impact on the performance of the SAEV fleet and its capacity to respond to these emergency situations. This can be directly implemented using the modeling and optimization approach proposed in this paper.  

\section{Conclusion}

In this paper, a dynamic optimization model for Shared-Autonomous Electric Vehicles operation for providing Autonomous Mobility-on-Demand services and Vehicle-to-Building (V2B) power services to critical buildings is proposed. The model considers the interdependence between the power demand of critical buildings and the mobility demands of passengers for an electrified mobility system. The model proposed explicitly considers the potential of a SAEV fleet to provide V2B power services to critical buildings, such as hospitals or other critical services, in cases of emergency with the aim of improving the building's energy resilience. To demonstrate the potential advantage of SAEV in improving the critical building energy resilience via V2B service, a case study based on the Ile-de-France region is considered, based on realistic  transport data and considering a hospital requiring emergency power within this region.

The results show that in normal conditions (i.e. no power outages in critical buildings) the SAEV fleet is capable of fully satisfying the passengers traveling demands with no significant waiting times. An appearance of an emergency power outage requiring V2B services may impact the mobility service quality as vehicles are contributing to the power demand at the expense of satisfying the mobility requirements. However, even if passenger total waiting times increase in these situations, the results show that adopting V2B services significantly improves the critical building’s resilience towards power outages, allowing it to recover at a much lower cost while maintaining acceptable levels of mobility service. Furthermore, it is shown that, for our case study, the SAEV fleet is capable of satisfying the emergency power service at a lower cost compared to investing in extra backup generators, unless the outage occurs at significantly high frequencies (more than 121 times a year). Finally, it is shown that considering fast charging and discharging can significantly improve the fleet  capability in providing high levels of mobility services in outage situations as a result of reducing the vehicle required to discharge at the same time. These results highlight that SAEV can, indeed, effectively be operated in such a way as to improve critical buildings' power resilience. Future studies can expand on these analyses to consider the inherent uncertainties in the system as well as the consideration of more detailed modeling features, such as traffic congestion and road capacities, among others. Particularly, our proposed dynamic modeling approach can be combined with adequate simulation models to generate random outages in critical building and investigate the expected performance of the system. In addition, further work could explore Reinforcement Learning algorithms as alternatives for solving the resulting sequential optimization problem as they have the potential of capturing elements of uncertainties more effectively, with fewer assumptions on the developed scenarios considered.

\appendix

\section{Results for vehicles trajectories of the SAEV fleet throughout one day under normal operating conditions}
\label{appen1}

\begin{figure}[h!]
\centering
    \subfloat[]{
		\includegraphics[width=0.42\textwidth]{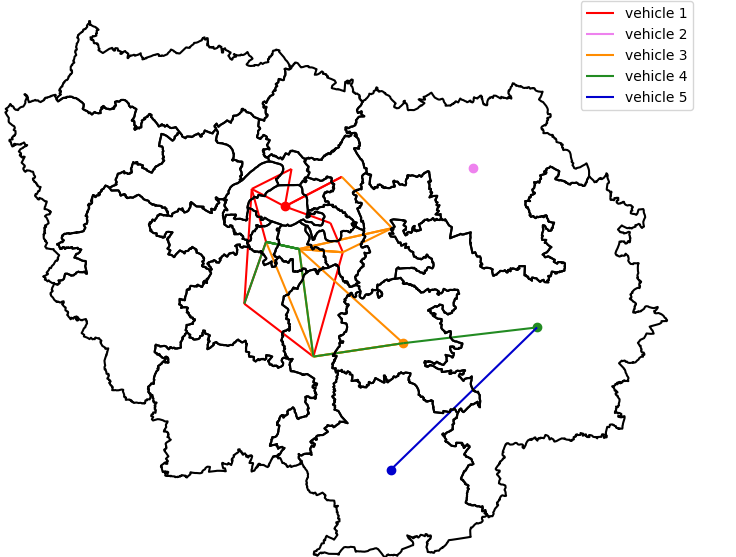}\hspace{-1pt}
		}
    \subfloat[]{
		\includegraphics[width=0.42\textwidth]{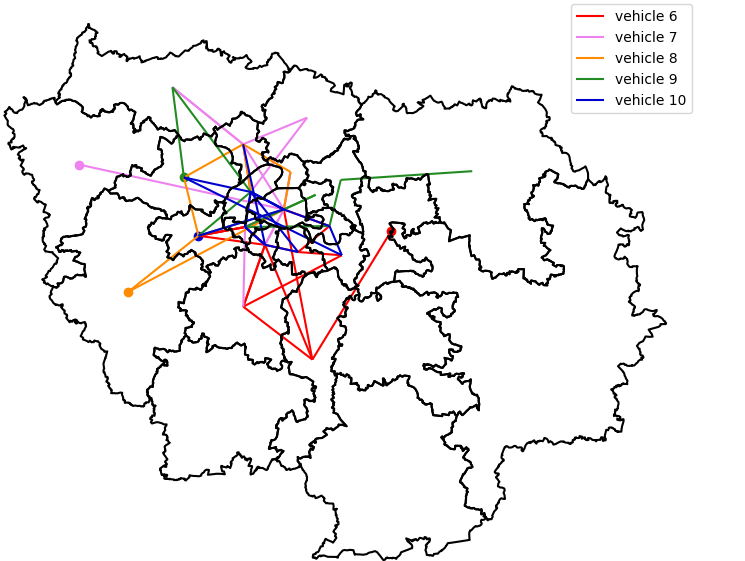}\vspace{2pt}
		}

\vfil
    \subfloat[]{
		\includegraphics[width=0.42\textwidth]{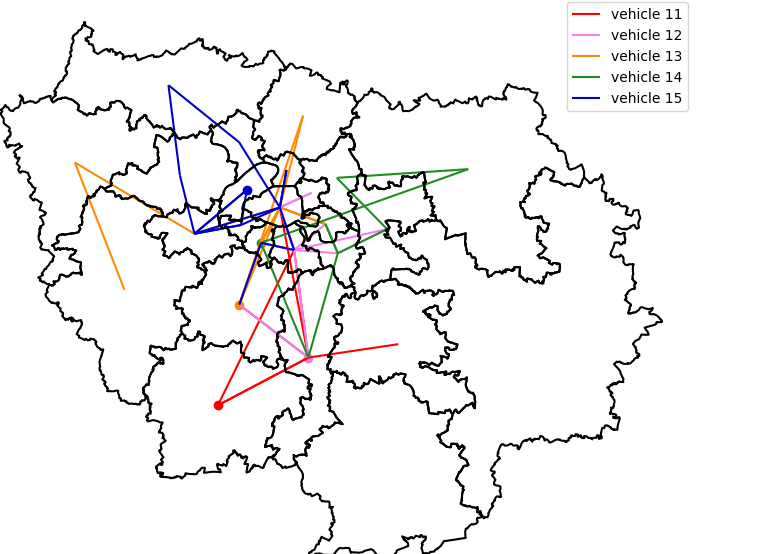}\hspace{-1pt}
		}
    \subfloat[]{
		\includegraphics[width=0.42\textwidth]{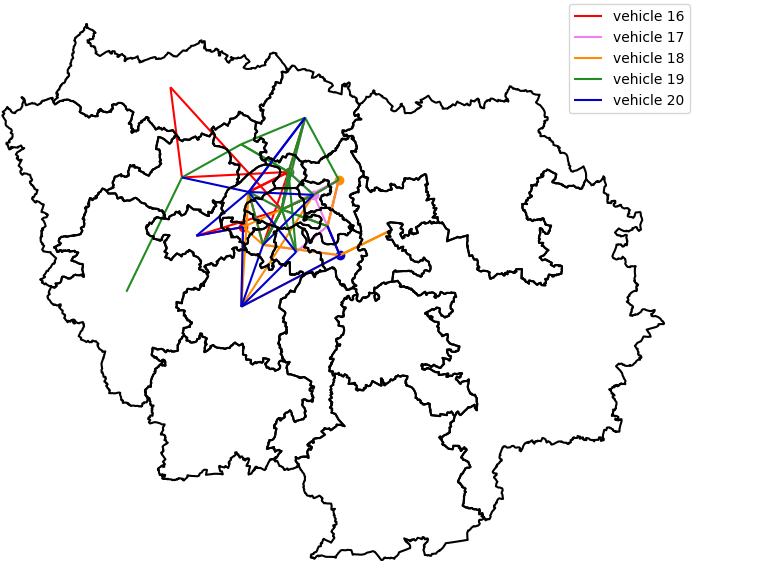}\vspace{2pt}
		}
\vfil
    \subfloat[]{
		\includegraphics[width=0.42\textwidth]{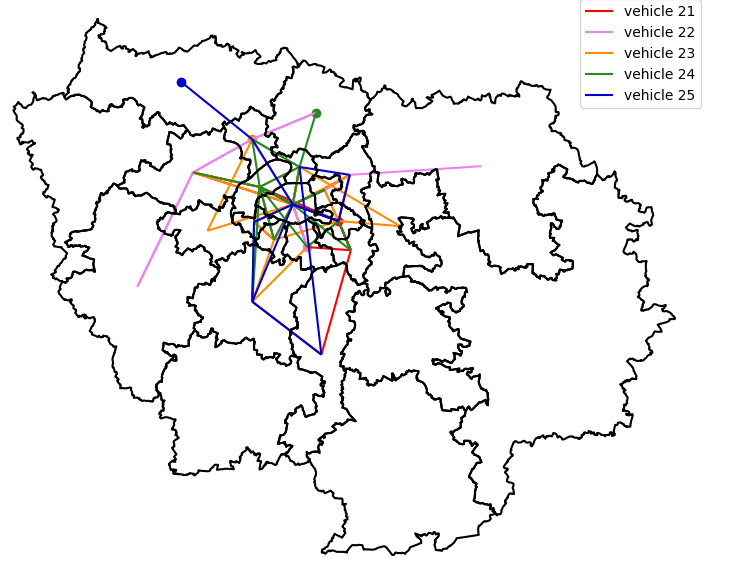}\hspace{-1pt}
		}
    \subfloat[]{
		\includegraphics[width=0.42\textwidth]{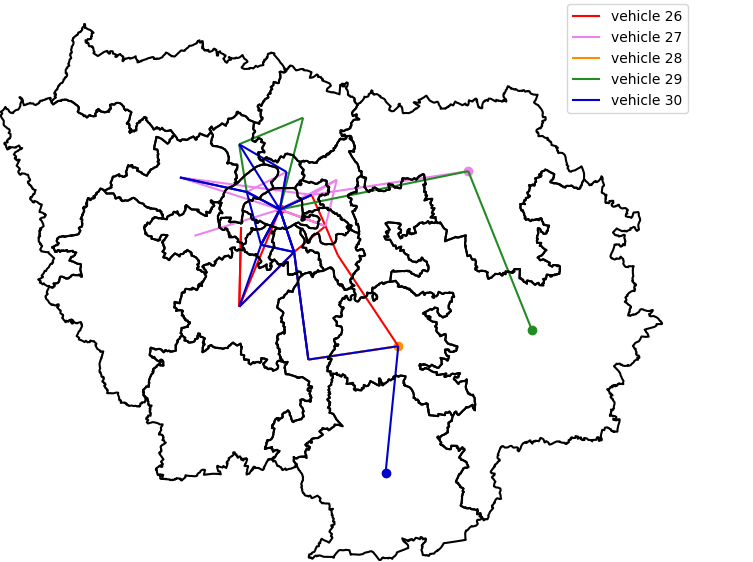}\vspace{2pt}
		}
    \caption{Results for vehicles trajectories throughout the day under normal operational conditions.}
    \label{appen:pic1}
\end{figure}

Figure~(\ref{appen:pic1}) illustrates the trajectories taken by all the vehicles of the SAEV fleet throughout the day as obtained by the optimal solution. Other than the general observation that most of the vehicles trajectories are concentrated within the central area of the Ile-de-France region, some vehicles (such as vehicle 2) may be idle in some scenarios where their services are not needed. They represent a backup for the system in these cases. Moreover, other vehicles (such as vehicle 5) may move only once during the day as a result of being far from the central area and only a few trip requests exist around.

\bibliographystyle{unsrtnat} 
{\small
\bibliography{bib}} 

\end{document}